\documentclass[]{article}
\usepackage{amsfonts}
\usepackage{tikz}
\usepackage{color}
\usepackage{amsmath}
\usepackage{MnSymbol}
\usepackage{standalone}
\usepackage{stmaryrd}
\usepackage{amsmath} 
\usepackage{pgfplots}
\usepackage{kbordermatrix}
\usepackage{pdfpages}
\usepackage{float}
\usepackage{geometry}
\usetikzlibrary{arrows}
\usetikzlibrary{calc}
\usetikzlibrary{matrix}
\usetikzlibrary{fit}
\usetikzlibrary{shapes.arrows}

\floatstyle{plain}
\newfloat{SIMAXtable}{th}{myt}
\floatname{SIMAXtable}{{\sc\footnotesize Table}}

\def\ER{\mathbb{R}}
\def\E{\mathbb{E}}
\def\f{\ensuremath{\mathbf{f}}}
\def\e{\ensuremath{\mathbf{e}}}
\def\u{\mathbf{u}}
\def\y{\mathbf{y}}

\def\V{\ensuremath{\mathbf{V}}}
\def\I{\ensuremath{\mathbf{I}}}
\def\X{\ensuremath{\mathbf{X}}}
\def\H{\ensuremath{\mathbf{H}}}

\def\T{\ensuremath{\mathbf{T}}}
\def\A{\ensuremath{\mathbf{A}}}
\def\W{\ensuremath{\mathbf{W}}}
\def\U{\ensuremath{\mathbf{U}}}
\def\w{\ensuremath{\mathbf{w}}}
\def\v{\ensuremath{\mathbf{v}}}
\def\a{\ensuremath{\mathbf{a}}}
\def\b{\ensuremath{\mathbf{b}}}
\def\c{\ensuremath{\mathbf{c}}}
\def\h{\ensuremath{\mathbf{h}}}

\def\Q{\ensuremath{\mathbf{Q}}}
\def\D{\ensuremath{\mathbf{D}}}
\def\L{\ensuremath{\mathbf{L}}}

\def\Ds{\ensuremath{\mathbf{D}_\Sigma}}
\def\Vs{\ensuremath{\mathbf{U}_\Sigma}}

\def\A{\ensuremath{\mathbf{A}}}

\def\g{\ensuremath{\mathbf{g}}}
\def\x{\ensuremath{\mathbf{x}}}
\def\z{\ensuremath{\mathbf{z}}}
\def\Jcal{\ensuremath{\mathcal{J}}}
\def\J{\ensuremath{\mathbf{J}}}
\def\Tcal{\ensuremath{\mathcal{T}}}
\def\vec#1{\mathrm{vec}\mathord{\bigl(#1\bigr)}}

\def\norm#1{\ensuremath{\left\| #1 \right\|}}
\def\normf#1{\norm{#1}_F}
\def\cpdgen#1#2#3{\llbracket#1\mathord{,}#2\mathord{,}#3\rrbracket}
\def\CPD{{\sc cpd}}

\def\rank{\mathrm{rank}}

\def\weight{\ensuremath{\mathord{\mathbf{\Omega}}}}
\def\matlab{\textsc{Matlab}}

\def\B{\ensuremath{\mathbf{B}}}

\def\quadtext#1{\quad\textrm{#1}\quad}

\def\weightone{\ensuremath{\weight_{(1)}}}
\def\weighttwo{\ensuremath{\weight_{(2)}}}
\def\weightthree{\ensuremath{\weight_{(3)}}}

\def\P{\ensuremath{\mathbf{P}}}
\def\eqref#1{Equation~(\ref{#1})}

\def\rSigma{\ensuremath{\mathord{\overline{r}}}}
\def\matlab{\textsc{matlab}}
\def\fig#1{\textsc{Fig}.~\ref{#1}}

\def\covarM{\Sigma_\Jcal}
\def\var{\mathop{\textrm{var}}}
\def\cpdgen#1#2#3{\llbracket#1\mathord{,}#2\mathord{,}#3\rrbracket}

\newcounter{algorCounter}
\newenvironment{algor}[1]%
{ \noindent\rule{\textwidth}{0.2ex} \stepcounter{algorCounter}\textbf{Algorithm \thealgorCounter}. #1\par} %
{ \par\noindent\rule{\textwidth}{0.2ex}\vskip-0.1\baselineskip }

\def\cov{\textrm{cov}}

\title{Weighted tensor decomposition for approximate decoupling of multivariate polynomials\thanks{This work was supported in part by the Fund for Scientific Research (FWO-Vlaanderen), the Flemish Government (Methusalem), the Belgian Government through the Interuniversity Poles of Attraction (IAP VII) Program, the ERC advanced grant SNLSID under contract 320378, the ERC starting grant SLRA under contract 258581, and FWO project G028015N.}} 

\author{Gabriel Hollander, Philippe Dreesen, Mariya Ishteva, Johan Schoukens\thanks{Department  of  VUB-ELEC,  Vrije  Universiteit  Brussel  (VUB),  B-1050  Brussels  (gabriel.hollander@vub.ac.be, philippe.dreesen@vub.ac.be, mariya.ishteva@vub.ac.be, johan.schoukens@vub.ac.be).}}

\begin{document}
\maketitle

\begin{abstract}
Multivariate polynomials arise in many different disciplines. Representing such a polynomial as a vector of univariate polynomials can offer useful insight, as well as more intuitive understanding. For this, techniques based on tensor methods are known, but these have only been studied in the exact case. In this paper, we generalize an existing method to the noisy case, by introducing a weight factor in the tensor decomposition. Finally, we apply the proposed weighted decoupling algorithm in the domain of system identification, and observe smaller model errors.
\end{abstract}

\pagestyle{myheadings}
\thispagestyle{plain}
\markboth{}{Weighted tensor decompositions}

\section{Introduction and notations}\label{sec:intro}
The starting point in this paper is a multivariate vector function $\f\colon\ER^m\to\ER^n$, where $f_i$ ($1 \le i \le n$) is a polynomial in $m$ variables of degree at most~$d$. The variables of~\f\ will be denoted as $\u = (u_1,\ldots,u_m)$ and the values as $\f(\u) = \y = (y_1,\ldots,y_n) = \bigl(f_1(\u),\ldots,f_n(\u)\bigr)$. This function may contain cross terms of monomials, for example $c_1u_1u_2$ or $c_2u_3^2u_2^2$, where $c_1,c_2\in\ER$, in which case it is called \emph{coupled}.

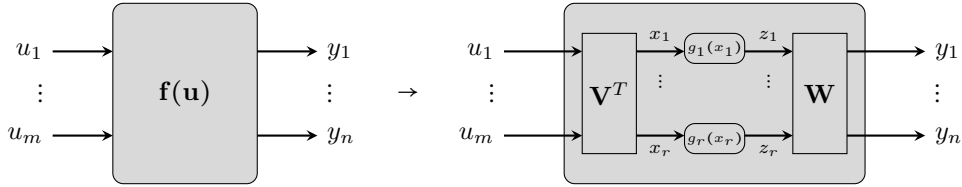
\begin{figure}[ht]
	\def\w{2.5}
\def\ww{1.2}
\def\wi{0.5}
\def\fign{1}
\def\u{\mathbf{u}}
\def\x{\mathbf{x}}
\def\J{\mathbf{J}}
\begin{tikzpicture}[>=stealth, scale=0.8]
    \filldraw[fill=gray!30, rounded corners = 4] (-\ww,-1.5) rectangle (\ww,1.5);
    \draw[->, thick] (-\ww-1, 0.7) node[anchor=east] {$u_1$} -- (-\ww, 0.7);
    \draw (-\ww-1, 0.0) node[anchor= east] {$\vdots$};
    \draw[->, thick] (-\ww-1, -0.7) node[anchor=east] {$u_m$} -- (-\ww, -0.7);
    \draw[->, thick] (\ww, 0.7) -- (\ww+1, 0.7) node[anchor=west] {$y_1$};
    \draw (\ww+1, 0.0) node[anchor= west] {$\vdots$};
    \draw[->, thick] (\ww, -0.7) -- (\ww+1, -0.7) node[anchor=west] {$y_n$};
    \draw (0,0) node {$\mathbf{f}(\u)$};
    \draw (\ww+2.5, 0) node {$\rightarrow$};
    \begin{scope}[xshift=8.8cm]
    \filldraw[fill=gray!30, rounded corners = 4] (-\w,-1.5) rectangle (\w,1.5);
    \draw[->, thick] (-\w-1, 0.7) node[anchor=east] {$u_1$} -- (-\w+0.3, 0.7);
    \draw (-\w-1, 0.0) node[anchor= east] {$\vdots$};
    \draw[->, thick] (-\w-1, -0.7) node[anchor=east] {$u_m$} -- (-\w+0.3, -0.7);
    \draw[->, thick] (\w-0.3, 0.7) -- (\w+1, 0.7) node[anchor=west] {$y_1$};
    \draw (\w+1, 0.0) node[anchor= west] {$\vdots$};
    \draw[->, thick] (\w-0.3, -0.7) -- (\w+1, -0.7) node[anchor=west] {$y_n$};
    \draw[rounded corners = 4] (-\wi,0.5) rectangle (\wi, 1);
    \draw[rounded corners = 4] (-\wi,-0.5) rectangle (\wi, -1);
    \draw (0, 0.75) node {$\scriptscriptstyle g_1(x_1)$};
    \draw (0, -0.75) node {$\scriptscriptstyle g_r(x_r)$};
    \draw (-\w+0.3,-1) rectangle (-\w+1.2,1);
    \draw (\w-1.2,-1) rectangle (\w-0.3,1);
    \draw (-\w+0.75,0) node {$\mathbf{V}^T$};
    \draw (\w-0.75,0) node {$\mathbf{W}$};
    \draw[->, thick] (-\w+1.2, 0.7) -- node[above] {$\scriptstyle x_1$} (-\wi, 0.7);
    \draw[->, opacity=0] (-\w+1.2, 0.2) -- node[opacity=1] {$\scriptstyle\vdots$} (-\wi, 0.2);
    \draw[->, thick] (-\w+1.2, -0.7) -- node[below] {$\scriptstyle x_r$} (-\wi, -0.7);
    \draw[<-, thick] (\w-1.2, 0.7) -- node[above] {$\scriptstyle z_1$} (\wi, 0.7);
    \draw[->, opacity=0] (\w-1.2, 0.2) -- node[opacity=1] {$\scriptstyle\vdots$} (\wi, 0.2);
    \draw[<-, thick] (\w-1.2, -0.7) -- node[below] {$\scriptstyle z_r$} (\wi, -0.7);
    \end{scope}
\end{tikzpicture}
	\caption{The decoupling process: given \f, find the matrices \V\ and \W\ and the univariate functions~$g_1, \ldots, g_r$.}
	\label{fig:dec}
\end{figure}

The principal goal of this article is to find a \emph{decoupled representation} of~\f\ as illustrated in \fig{fig:dec}.
Given $\f$, we wish to find transformation matrices $\V\in\ER^{m\times r}$ and~$\W\in\ER^{n\times r}$ and a vector~$\g(\x)=\bigl(g_1(x_1), \ldots, g_r(x_r)\bigr)$ of univariate polynomials, such that,
$$
\f(\u) \approx \W\g(\V^T\u),
$$
for all inputs $\u$. The first internal variable is denoted as $\x=\V^T\u$, and the second as $\z = \g(\V^T\u)$. Furthermore, the number~$r$ of internal \emph{branches} is assumed to be predefined. This decoupled representation offers a way to study~\f\ without cross terms, which can be an advantage for certain applications, as it helps its physical or intuitive understanding.

\begin{figure}[ht]
\pgfmathsetmacro{\rr}{1.4}
\pgfmathsetmacro{\nn}{1}
\pgfmathsetmacro{\mm}{0.6}
\pgfmathsetmacro{\NN}{2.5}
\pgfmathsetmacro{\wVector}{0.25}
\pgfmathsetmacro{\sep}{0.3}
\pgfmathsetmacro{\vsep}{0.15}
\pgfmathsetmacro{\dnshift}{-3cm}
\pgfmathsetmacro{\horshift}{-2cm}
\centerline{\begin{tikzpicture}[scale=0.7,x={(1cm,0cm)},y={(0cm,1cm)},z={(-3.85mm,-3.85mm)},thick,xshift=-0.5cm]
\begin{scope}[xshift=1.5*\horshift]
\draw (0,0.6*\rr,0) node {$=$};
\end{scope}
\begin{scope}[xshift=1.5*\horshift,yshift=\dnshift]
\draw (0,0.6*\rr,0) node {$=$};
\end{scope}
\begin{scope}[xshift=3*\horshift]
    \begin{scope}
    \draw [line width=2pt,draw=red,fill=black!45!red] ($(0,0,-4*2*\sep)$) -- ++(\mm,0,0) -- ++(0,\nn,0) --++ (-\mm,0,0) -- cycle;
    \draw [line width=2pt,draw=green,fill=black!45!green] ($(0,0,-3*2*\sep)$) -- ++(\mm,0,0) -- ++(0,\nn,0) --++ (-\mm,0,0) -- cycle;
    \draw [line width=2pt,draw=cyan,fill=black!45!cyan] ($(0,0,-2*2*\sep)$) -- ++(\mm,0,0) -- ++(0,\nn,0) --++ (-\mm,0,0) -- cycle;
    \draw [line width=2pt,draw=magenta,fill=black!45!magenta] ($(0,0,-1*2*\sep)$) -- ++(\mm,0,0) -- ++(0,\nn,0) --++ (-\mm,0,0) -- cycle;
    \draw [line width=2pt,draw=blue,fill=black!45!blue] (0,0,0) -- ++(\mm,0,0) -- ++(0,\nn,0) --++ (-\mm,0,0) -- cycle;
    \end{scope}
\end{scope}
\begin{scope} 
    \begin{scope}
    \draw [draw=black,fill=black!15] (-\sep,\rr,0) -- ++(-\rr,0,0) -- ++(0,-\nn,0) -- ++(\rr,0,0) -- cycle;
    \draw (-\sep-0.5*\rr,\rr-0.5*\nn,0) node {$\mathbf{W}$}; 
    \end{scope}
    \begin{scope}
    \draw [draw=white,fill=white,fill opacity=0.7] ($(0,0,-4*2*\sep)$) -- ++(\rr,0,0) -- ++(0,\rr,0) --++ (-\rr,0,0) -- cycle;
    \draw [draw=red,line width=3pt,triangle 90 cap-triangle 90 cap,color=red] ($(0,\rr,-4*2*\sep)$) -- ++(\rr,-\rr,0);
    \draw [draw=black] ($(0,0,-4*2*\sep)$) -- ++(\rr,0,0) -- ++(0,\rr,0) --++ (-\rr,0,0) -- cycle;
    \draw [draw=white,fill=white,fill opacity=0.7] ($(0,0,-3*2*\sep)$) -- ++(\rr,0,0) -- ++(0,\rr,0) --++ (-\rr,0,0) -- cycle;
    \draw [draw=green,line width=3pt,triangle 90 cap-triangle 90 cap,color=green] ($(0,\rr,-3*2*\sep)$) -- ++(\rr,-\rr,0);
    \draw [draw=black] ($(0,0,-3*2*\sep)$) -- ++(\rr,0,0) -- ++(0,\rr,0) --++ (-\rr,0,0) -- cycle;
    \draw [draw=white,fill=white,fill opacity=0.7] ($(0,0,-2*2*\sep)$) -- ++(\rr,0,0) -- ++(0,\rr,0) --++ (-\rr,0,0) -- cycle;
    \draw [draw=cyan,line width=3pt,triangle 90 cap-triangle 90 cap,color=cyan] ($(0,\rr,-2*2*\sep)$) -- ++(\rr,-\rr,0);
    \draw [draw=black] ($(0,0,-2*2*\sep)$) -- ++(\rr,0,0) -- ++(0,\rr,0) --++ (-\rr,0,0) -- cycle;
    \draw [draw=white,fill=white,fill opacity=0.7] ($(0,0,-1*2*\sep)$) -- ++(\rr,0,0) -- ++(0,\rr,0) --++ (-\rr,0,0) -- cycle;
    \draw [draw=magenta,line width=3pt,triangle 90 cap-triangle 90 cap,color=magenta] ($(0,\rr,-1*2*\sep)$) -- ++(\rr,-\rr,0);
    \draw [draw=black] ($(0,0,-1*2*\sep)$) -- ++(\rr,0,0) -- ++(0,\rr,0) --++ (-\rr,0,0) -- cycle;
    \draw [draw=white,fill=white,fill opacity=0.7] ($(0,0,-0*2*\sep)$) -- ++(\rr,0,0) -- ++(0,\rr,0) --++ (-\rr,0,0) -- cycle;
    \draw [draw=blue,line width=3pt,triangle 90 cap-triangle 90 cap,color=blue] ($(0,\rr,-0*2*\sep)$) -- ++(\rr,-\rr,0);
    \draw [draw=black] ($(0,0,-0*2*\sep)$) -- ++(\rr,0,0) -- ++(0,\rr,0) --++ (-\rr,0,0) -- cycle;
    \end{scope}
    \begin{scope}
        \coordinate (startingpointrightfactor) at ($(\rr +0.385*4*2*\sep + \sep,\rr,0)$); 
     (0,0,0) -- ++(\rr,0,0) -- ++(0,0,-\rr) -- ++(0,\rr,0) -- ++(-\rr,0,0) -- ++(0,0,\rr) -- cycle;
    \draw[draw=black,fill=black!15] (startingpointrightfactor) -- ++ (\rr,0,0) -- ++(0,-\mm,0) -- ++(-\rr,0,0) -- cycle;
    \draw ($(startingpointrightfactor) + (0.5*\rr,-0.5*\mm,0)$) node {$\mathbf{V}$}; 
    \end{scope}
\end{scope}
\begin{scope}[xshift=-1cm,yshift=0.8*\dnshift]
    \begin{scope}
        \draw (-4.7,2) node {$\mathcal{J}$};
        \draw[draw=black,fill=black!15] (0,-\vsep,0) -- ++ (0,-\nn,0) -- ++ (\wVector,0,0) -- ++ (0,\nn,0) --  cycle;
        \draw (-2*\wVector,-0.8*\vsep-0.8*\nn,0) node {$\mathbf{w}_1$};
        \draw[draw=black,fill=black!15] (\wVector+\vsep,0,0) -- ++ (\mm,0,0) -- ++ (0,\wVector,0) -- ++ (-\mm,0,0) -- cycle;
        \draw (\wVector+\vsep+0.8*\mm,-2*\wVector,0) node {$\mathbf{v}_1$};
        \draw[draw=black,fill=black!15] (0,\wVector,-\vsep) -- ++ (0,0,-\NN) -- ++ (\wVector,0,0) -- ++ (0,0,\NN) -- cycle;
        \draw (-3*\wVector,\wVector,-0.9*\NN) node {$\mathbf{h}_1$};
        \coordinate (pluscoord1) at (0.385*\NN+3*\sep,-0.2*\nn,0); 
        \draw (pluscoord1) node {$+$};
        \coordinate (startcoord2) at ($(pluscoord1) + (3*\sep,0,0)$);
        \draw (startcoord2) node {$\ldots$};
        \coordinate (startcoordr) at ($(startcoord2) +  (3*\sep,0,0)$);
        \draw (startcoordr) node {$+$};
        \coordinate (startcoordlast) at ($(startcoordr)+  (3*\sep,0,0)$);
        \draw[draw=black,fill=black!15] ($(startcoordlast) + (0,-\vsep,0)$) -- ++ (0,-\nn,0) -- ++ (\wVector,0,0) -- ++ (0,\nn,0) --  cycle;
        \draw ($(startcoordlast) + (-2*\wVector,-0.8*\vsep-0.8*\nn,0)$) node {$\mathbf{w}_r$};
        \draw[draw=black,fill=black!15] ($(startcoordlast) + (\wVector+\vsep,0,0)$) -- ++ (\mm,0,0) -- ++ (0,\wVector,0) -- ++ (-\mm,0,0) -- cycle;
        \draw ($(startcoordlast) + (\wVector+\vsep+0.8*\mm,-2*\wVector,0)$) node {$\mathbf{v}_r$};
        \draw[draw=black,fill=black!15]  ($(startcoordlast) + (0,\wVector,-\vsep)$) -- ++ (0,0,-\NN) -- ++ (\wVector,0,0) -- ++ (0,0,\NN) -- cycle;
        \draw ($(startcoordlast) + (-3*\wVector,\wVector,-0.9*\NN)$) node {$\mathbf{h}_r$};
    \end{scope}
\end{scope}
\end{tikzpicture}}
\caption{Graphical representation of the core of Algorithm~1.}
\label{fig_decouplingGraphical}\end{figure}
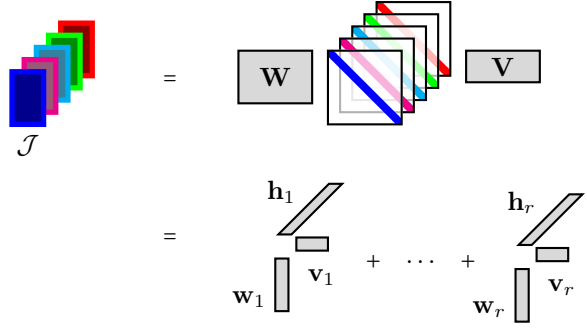

To our knowledge, \cite{Dreesen2015} and~\cite{MaartenSchoukens2014} offer a solution to this problem under the special assumption that an exact decomposition with~$r$ branches exists. Under the extra condition of homogeneous polynomials, this has also been studied in~\cite{Tiels2013}. Section~1.2 of \cite{Dreesen2015} also refers to the related Waring problem. In the case of state-space models, this problem is addressed in~\cite{Mulders2014}. Because the solution of~\cite{Dreesen2015} seems to be computationally easier, we have chosen to use and generalize this algorithm, which is based on the first-order derivative information of \f. At its core, tensor decompositions are used and the method is outlined in Algorithm~1. \fig{fig_decouplingGraphical} shows a graphical representation. An overview of tensor decompositions can be found in~\cite{Comon2002}, \cite{Comon2009a}, \cite{Kolda2009} and \cite{Lathauwer2000}.

\begin{algor}{Decomposing a multivariate polynomial \f\ having an exact decomposition. In this section, we shortly introduce the algorithm of \cite{Dreesen2015}.}
\begin{enumerate}
	\item Evaluate the Jacobian matrix of \f
	$$
	\J(\u) = \begin{bmatrix}
	\frac{\partial f_1}{\partial u_1}(\u) & \cdots & \frac{\partial f_1}{\partial u_m}(\u)\\
	\vdots & \ddots & \vdots \\
	\frac{\partial f_n}{\partial u_1}(\u) & \cdots & \frac{\partial f_n}{\partial u_m}(\u)
	\end{bmatrix}
	$$
	in $N$ randomly chosen points $\u^{(1)},\ldots,\u^{(N)}$. The number $N$ of sampling points is chosen by the user. The equality $\f(\u) = \W\g(\V^T\u)$ implies for the Jacobians that
	$$
		\J(\u) = \W
	\begin{bmatrix}
		g_1'(\v_1^T\u) & & 0\\
		& \ddots & \\
		0 & & g_r'(\v_r^T\u) 
	\end{bmatrix}\V^T,
	$$
	where the matrix of derivatives is zero outside of the diagonal and $\v_1,\ldots,\v_r$ denote the columns of $\V$, and $g'_i$ denote the derivative of the $i$-th component of $\g$.
	\item Stack the Jacobians into a three-way tensor $\Jcal$ of dimensions $n\times m\times N$. Here, the $k$-th frontal slice of \Jcal\ consists of the first-order information of \f\ evaluated in the sampling point $\u^{(k)}$, i.e., $\Jcal(\colon,\colon,k) = \J(\u^{(k)})$ (where $1\le k\le N$). Here, the \matlab-notation $\Jcal(\colon,\colon,k)$ is used for the $k$-th frontal slide of \Jcal.
	\item Compute the Canonical Polyadic Decomposition (CPD) of
	\begin{equation}
		\Jcal = \sum_{i=1}^r \w_i\circ\v_i\circ\h_i.
	\label{eq_cpd1}\end{equation}
	Here, the vectors $\w_i$ (respectively $\v_i$) define the columns of the matrix~\W\ (respectively \V) of the decoupled representation. Furthermore, the vectors $\h_1,\ldots,\h_r$ contain the first-order information of the internal univariate functions $g_1(x_1),\ldots,g_r(x_r)$ evaluated in the $N$ sampling points, after transformation by $\V^T$, i.e., $\x^{(k)} = \V^T\u^{(k)}$. We thus have $h_{ij} = g'_j(\v^T_j\u^{(i)})$.
	\item Starting from the vectors $\h_1,\ldots,\h_r$, reconstruct the internal univariate functions $g_1(x_1),\ldots,g_r(x_r)$. This works by fitting the derivatives $g'_j$ using the vectors $\h_1,\ldots,\h_r$, and then to recover the functions $g_j$ with an integration step. This method is described in Section~II.C of~\cite{Dreesen2015a}. Since we focus our attention to noisy coefficients of \f, it seems reasonable to use all the information available in the vectors $\h_1,\ldots,\h_r$, instead of the method proposed in Section~2.4 of~\cite{Dreesen2015}.
\end{enumerate}
\end{algor}

In~\cite{Dreesen2015}, it is shown that Algorithm~1 works well in case that the function~\f\ admits an exact decomposition. In this paper, we generalize the method to the noisy case: here, it is not assumed that an exact decoupling of \f\ exists, and instead, we will search for an approximated decoupling. In this regard, the coefficients of \f\ are thus considered ``noisy'': $\f = \f_0 + \v_\f$, where $\f_0$ has an exact decoupling and $\v_\f$ is zero-mean noise. 

In order to decouple this noisy coupled function \f, the covariance matrix~$\mathbf{\Sigma}_\f = \cov(\v)$ of~\f\ will be assumed to be known throughout this paper. Because this matrix is easy to approximate when doing numerical experiments or measurements, this seems a reasonable assumption. The matrix~$\mathbf{\Sigma}_\f$ contains the variances of and covariances between the different coefficients. This will lead to the creation of a weight matrix to be used during the decoupling process. This weight matrix will be defined, as common in a weighted least squares approximation problem, as the (pseudo) inverse of a covariance matrix. In conclusion, returning to the outline of the decoupling method, the attention in this paper will be focused on step~3 of Algorithm~1: the CPD will be generalized to a \emph{weighted CPD}.

This paper is organized as follows: in Section~\ref{sec:weight}, the ``covariance'' matrix~$\mathbf{\Sigma}_\Jcal$ of the Jacobian elements will be constructed as a linear transformation of the matrix $\mathbf{\Sigma}_\f$. In Section~\ref{sec:wCPD}, the generalization of the CPD with weights will be discussed. Finally, Section~\ref{sec:results} summarizes the numerical experiments of the weighted CPD and shows an application of the results in the domain of system identification, while Section~\ref{sec:conclusion} contains the conclusion and ideas for future work.

\section{Constructing the covariance matrix~$\mathbf{\Sigma}_\Jcal$}\label{sec:weight}
In order to create a weighted CPD, the covariance matrix~$\mathbf{\Sigma}_\Jcal$ of the Jacobian elements will be constructed as a linear transformation of~$\mathbf{\Sigma}_\f$, the covariance matrix of the coefficients of \f, except the constant terms. Because \Jcal\ contains $mnN$ elements, $\mathbf{\Sigma}_\Jcal$ will have dimensions $(nmN)\times(nmN)$.

In practice, three different matrices (and hence, three different weight matrices) will be discussed in Section~\ref{sec:wCPD}: (1) only the variances of each Jacobian element will be taken into account while the other covariances will be set to~0, (2) the slice-wise covariances will be taken into account as well, keeping the rest as 0, and, (3) all the covariances of $\Jcal$ will be used, forming a dense covariance matrix. This way, even though the element-wise and slice-wise defined matrices are approximations of the full covariance matrix, and are not by themselves well-defined covariance matrices, we will still use this term to denote them. Furthermore, we will denote them respectively $\mathbf{\Sigma}_\Jcal^\textrm{e}$, $\mathbf{\Sigma}_\Jcal^\textrm{s}$ and $\mathbf{\Sigma}_\Jcal^\textrm{d}$, and will use $\mathbf{\Sigma}_\Jcal$ if we wish to denote any one of them.

For ease of reading, we will often illustrate dimensions and values for the special case where $m=n=d=2$. The number $\ell$ of monomials of degree at most $d$, given by $\binom{m+d}{m}$, is in this case 6, and the coupled function \f\ can be written as
$$
\f(u_1,u_2) = 
\begin{bmatrix}
f_1(u_1,u_2) \\
f_2(u_1,u_2) 
\end{bmatrix}
=
\begin{bmatrix}
   	c_1 + c_2 u_1 + c_3 u_2 + c_4 u_1^2 + c_5 u_1u_2 + c_6 u_2^2 \\
   	d_1 + d_2 u_1 + d_3 u_2 + d_4 u_1^2 + d_5 u_1u_2 + d_6 u_2^2
\end{bmatrix},
$$
where we use $c_i$ and $d_i$ ($1\le i\le6$) for the coefficients of the first and second output of~\f, respectively. The covariance matrix $\mathbf{\Sigma}_\f$ of these coefficients is assumed to be known, and this is a $10\times10$ matrix. It is defined as follows:
\renewcommand{\kbldelim}{[}
\renewcommand{\kbrdelim}{]}
$$
\mathbf{\Sigma}_\f = \kbordermatrix{
    	& c_2 & \cdots & c_6 & d_2 & \cdots & d_6 \\
    c_2 & \var(c_2) & \cdots & \cov(c_2,c_6) & \cov(c_2,d_2) & \cdots & \cov(c_2,d_6)\\
    \vdots & \vdots & \ddots & \vdots & \vdots & \ddots & \vdots \\
    c_6 & \cov(c_6,c_2) & \cdots & \var(c_6) & \cov(c_6,d_2) & \cdots & \cov(c_6,d_6)\\
    d_2 & \cov(d_2,c_2) & \cdots & \cov(d_2,c_6) & \var(d_2) & \cdots & \cov(d_2,d_6)\\
    \vdots & \vdots & \ddots & \vdots & \vdots & \ddots & \vdots \\
    d_6 &  \cov(d_6,c_2) & \cdots & \cov(d_6,c_6) & \cov(d_6,d_2) & \cdots & \var(d_6) \\
  }.
$$
Because the first-order information of \f\ will be used, the information about the constant terms is not included in the matrix $\mathbf{\Sigma}_\f$. In general, the dimensions of $\mathbf{\Sigma}_\f$ are $\bigl((\ell-1) n\bigr)\times\bigl((\ell-1) n\bigr)$.

\subsection{Constructing the element-wise covariance matrix $\mathbf{\Sigma}_\Jcal^\textrm{e}$}
Using the notations of the previous section, it is possible to write the first Jacobian element evaluated at a point $\u^{(k)}$ ($1\le k\le N$) as follows:
$$
\J_{11}(\u^{(k)}) = \frac{\partial f_1}{\partial u_1}(\u^{(k)}) 
=
\begin{bmatrix}
1\; 0\; 2u_1^{(k)}\; u_2^{(k)}\; 0\
\end{bmatrix} 
\begin{bmatrix}
c_2 \\
\vdots \\
c_{6}
\end{bmatrix}.
$$
Derivatives of polynomials do not depend of the constant term $c_1$. It then follows for the variance of the element $\J_{11}(\u^{(k)})$ of \Jcal\ that
\begin{equation}
\var\bigl(\J_{11}(\u^{(k)})\bigr) = \begin{bmatrix}
1\; 0\; 2u_1^{(k)}\; u_2^{(k)}\; 0\;
\end{bmatrix} \mathbf{\Sigma}_{c} 
\begin{bmatrix}
1 \\ 
0 \\
2u_1^{(k)} \\
u_2^{(k)} \\
0
\end{bmatrix},
\label{eq_elementwise}\end{equation}
where
$$
\mathbf{\Sigma}_{c} = 
 \begin{bmatrix}
	 \var(c_2) & \cdots & \cov(c_2,c_6) \\
    \vdots & \ddots & \vdots \\
     \cov(c_6,c_2) & \cdots & \var(c_6) & \\
  \end{bmatrix}
$$
is the portion of $\mathbf{\Sigma}_\f$ consisting of the variances of and covariances between the coefficients~$c_2, \ldots, c_6$. When repeating this in a similar way for the other Jacobian elements, one finds all the variances of the elements of \Jcal, which, once collected into a matrix, give the element-wise matrix $\mathbf{\Sigma}_\Jcal^\textrm{e}$:
\begin{equation}
\mathbf{\Sigma}_\Jcal^\textrm{e} = \begin{bmatrix}
\var\J_{11}^{(1)} \\
& \var\J_{21}^{(1)} \\
&& \var\J_{12}^{(1)} \\
&&& \var\J_{22}^{(1)} \\
&&&& \ddots \\
&&&&&\var\J_{11}^{(N)} \\
&&&&&& \var\J_{21}^{(N)} \\
&&&&&&& \var\J_{12}^{(N)} \\
&&&&&&&& \var\J_{22}^{(N)} \\
\end{bmatrix}.
\label{eq_elementwiseCovar}\end{equation}
Here, $\J_{11}^{(k)}$ is used as shorthand notation for $\J_{11}(\u^{(k)})$, and all empty spaces in the matrix denote zero. Hence, the covariance matrix obtained when considering solely at the variances of Jacobian elements is diagonal. Implementing the weighted CPD decomposition using an element-wise weight has been described in~\cite{Andersson2000} or in \cite{Bader2015} (in the case of incomplete data). In the following two sections, this will be generalized to covariance matrices with more cross-covariances taken into account.

\subsection{Constructing the slice-wise covariance matrix $\mathbf{\Sigma}_\Jcal^\textrm{s}$}
As a next step, the covariances between Jacobian elements evaluated at a single sampling point $\u^{(k)}$ are taken into account for the construction of the matrix $\covarM$, but the covariances over several $\u^{(k)}$ are neglected. This is done by noting that the Jacobian elements of sampling point~$\u^{(k)}$ can be written as linear combinations of the monomials of degree at least one:
\setlength{\arrayrulewidth}{1pt}
$$
\begin{bmatrix}
\J_{11}(\u^{(k)})\\ \J_{21}(\u^{(k)}) \\
\J_{12}(\u^{(k)}) \\ \J_{22}(\u^{(k)})
\end{bmatrix} = 
\left[\begin{array}{ccccc| ccccc}
1 & 0 & 2u_1^{(k)} & u_2^{(k)} & 0 & 0 & 0 & 0 & 0 & 0 \\
0 & 0 & 0 & 0 & 0 & 1 & 0 & 2u_1^{(k)} & u_2^{(k)} & 0 \\
0 & 1 & 0 & u_1^{(k)} & 2u_2^{(k)} & 0 & 0 & 0 & 0 & 0 \\
0 & 0 & 0 & 0 & 0 & 0 & 1 & 0 & u_1^{(k)} & 2u_2^{(k)} \\
\end{array}\right]
\left[\begin{array}{@{}c@{}}
c_2 \\
\vdots \\
c_6 \\
\hline
d_2 \\
\vdots \\
d_6
\end{array}\right].
$$
This can be denoted shortly as 
\begin{equation}
\vec{\J(\u^{(k)})} = \A(\u^{(k)}) \left[\begin{array}{@{}c@{}}
c_2 \\
\vdots \\
c_6 \\
\hline
d_2 \\
\vdots \\
d_6
\end{array}\right],
\label{eq:defAmatrix}
\end{equation}
where, the matrix $\A\in\ER^{4\times10}$ (in general, $\A\in\ER^{(m n)\times((\ell-1) n)}$) contains the linear relationships between coefficients of the monomials and Jacobian elements. In the case that $d>2$, the matrix $\A$ may of course contain higher powers of the $u_i^{(k)}$. It follows that the covariance matrix of the Jacobian elements of sampling point~$\u^{(k)}$ is given by $\A(\u^{(k)}) \mathbf{\Sigma}_\f \A(\u^{(k)})^T \in \ER^{(mn)\times(mn)}$. When we repeat this for all the sampling points, we obtain the slice-wise covariance matrix~$\mathbf{\Sigma}_\Jcal^\textrm{s}$
\begin{equation}
\mathbf{\Sigma}_\Jcal^\textrm{s} = \begin{bmatrix}
\A(\u^{(1)}) \mathbf{\Sigma}_\f  \A(\u^{(1)})^T \\
& \A(\u^{(2)}) \mathbf{\Sigma}_\f \A(\u^{(2)})^T \\
&& \ddots \\
&&& \A(\u^{(N)}) \mathbf{\Sigma}_\f \A(\u^{(N)})^T
\end{bmatrix},
\label{eq_slicewise}\end{equation}
which is a block-diagonal matrix. Once again, the empty spaces in this matrix denote zero.

\subsection{Constructing the dense covariance matrix $\mathbf{\Sigma}_\Jcal^\textrm{d}$}
Finally, we generalize the slice-wise covariance matrix to the dense case, considering the (co)variances in all the elements of \Jcal. With the same definition for the matrix \A\ as in \eqref{eq:defAmatrix}, the dense covariance matrix~$\mathbf{\Sigma}_\Jcal^\textrm{d}$ is defined as
\begin{equation}
	\mathbf{\Sigma}_\Jcal^\textrm{d} = \begin{bmatrix}
	\A(\u^{(1)}) \\
	\vdots \\
	\A(\u^{(N)})
	\end{bmatrix} \mathbf{\Sigma}_\f\begin{bmatrix} \A(\u^{(1)})^T \cdots \A(\u^{(N)})^T \end{bmatrix}.
\label{eq_densecovar}\end{equation}
It follows immediately from this definition that, since~$\mathbf{\Sigma}_\Jcal^\textrm{d}\in\ER^{(mnN)\times(mnN)}$ and $\mathbf{\Sigma}_\f\in\ER^{((\ell-1)n)\times((\ell-1)n)}$, $\mathbf{\Sigma}_\Jcal^\textrm{d}$ is rank-deficient whenever $N>\frac{\ell-1}{m}$: its rank is then bounded by $(\ell-1)n$. This will lead to two different weighted CPD decompositions. On the one hand, decompositions using element-wise or slice-wise weights will be discussed in Section~\ref{sec:permutation} because the covariance matrix (and hence the weight matrix) has full rank; on the other hand, the dense-weight decompositions with a rank-deficient weight matrix will be discussed in Section~\ref{sec:f}.

\section{Computing the weighted CPD}\label{sec:wCPD}
Given the $n\times m\times N$-tensor \Jcal\ of Jacobian elements, it is possible to (approximately) write it as a sum of~$r$ rank-one tensors
\begin{equation}
	\Jcal \approx \sum_{i=1}^r \w_i\circ\v_i\circ\h_i.
\label{eq:sumRank1}\end{equation}
Here, the notation $\circ$ is used to denote the so-called outer product and it is defined as follows. Given vectors $\a = (a_1, \ldots, a_{n_a})$, $\b = (b_1, \ldots, b_{n_b})$ and $\c = (c_1, \ldots, c_{n_c})$, the outer product $\a\circ\b\circ\c$ is the $n_a\times n_b\times n_c$ tensor whose element in position $(i,j,k)$ ($1\le i\le n_a$, $1\le j\le n_b$, $1\le k\le n_c$) is given by
$$
(\a\circ\b\circ\c)_{i,j,k} = a_i b_j c_k.
$$
The outer product $\a\circ\b\circ\c$ is said to have \emph{rank one}. In expression~(\ref{eq:sumRank1}), the tensor~\Jcal\ is said to be approximated by a sum of rank~one tensors. We will use the same notation as in~\cite{Kolda2009} and will denote this sum as $\cpdgen\W\V\H$. 

Finding the matrices $\V$, $\W$ and $\H$ leads to optimizing the following non-linear cost function
$$
\begin{aligned}
	&\min_{\V,\W,\H} \normf{\Jcal - \cpdgen\W\V\H}^2 \\
	=& 
	\min_{\V,\W,\H} \bigl(\vec{\Jcal} - \vec{\cpdgen\W\V\H}\bigr)^T \bigl(\vec{\Jcal} - \vec{\cpdgen\W\V\H}\bigr),
\end{aligned}
$$
where $\vec\X$ denotes the vectorization of the matrix or tensor \X\ to a column vector. While this optimization problem is used for computing the unweighted CPD decomposition, the weighted CPD can be computed using a \emph{weighted norm} $\norm{\cdot}_{\weight}$ in the cost function, which includes a weight matrix $\weight$ as follows,
\begin{equation}
\begin{aligned}
	&\min_{\V,\W,\H} \norm{\vec{\Jcal} - \vec{\cpdgen\W\V\H}}^2_{\weight} \\
	=& \min_{\V,\W,\H} \left(\vec{\Jcal} - \vec{\cpdgen\W\V\H}\right)^{\!\!T} \weight \left(\vec{\Jcal} - \vec{\cpdgen\W\V\H}\right).
\end{aligned}
\label{eq:costfunction}\end{equation}
Optimizing this expression is described in the next sections, once for the element-wise and slice-wise weight, and once for the dense weight matrix. Naively, the weight matrix~$\weight$ can be intuitively seen as the (pseudo-) inverse of the covariance matrix~$\mathbf{\Sigma}_\Jcal$. This will be detailed in the following sections.
\label{sec:d_bd}
\subsection{Using element-wise and slice-wise weights}\label{sec:permutation}
In order to optimize the nonlinear cost function~(\ref{eq:costfunction}), we will use the workhorse method for computing the (unweighted) CPD: the Alternating Least Squares method described in~\cite{Carroll1970}, \cite{Harshman1970} and~\cite{Kolda2009}. This works by optimizing only one of the factors $\V$, $\W$ and $\H$ at a time, while keeping the two others fixed. Then, by alternating the optimized factor iteratively, a solution of the optimization problem can be found. Optimizing one factor (in the unweighted case) leads to three different least squares optimizations:
\begin{equation}
\min_{\W} \normf{\Jcal^T_{(1)} - (\H\odot\V) \W^T}^2,
\label{eq_als1}\end{equation}
$$
\min_{\V} \normf{\Jcal^T_{(2)} - (\H\odot\W) \V^T}^2,
$$
and
$$
\min_{\H} \normf{\Jcal^T_{(3)} - (\V\odot\W) \H^T}^2.
$$
Here, $\odot$ denotes the Khatri-Rao product (see, among others, \cite{Liu2008}), and $\Jcal_{(1)}, \Jcal_{(2)}, \Jcal_{(3)}$ denote the matricizations of the tensor \Jcal\ to the first, second and third mode, respectively. For this, we use the matricization order described in~\cite{Kolda2009}.

In order to take the weight matrix~$\weight$ into account, these (unweighted) least squares problems will be generalized to a set of weighted least squares optimizations, leading to a Weighted Alternating Least Squares method for finding the weighted CPD. Because this is done in a similar way for the three factors, more attention will be given to the update of $\W$ and only the differences with the two other factor updates will be emphasized.

\subsubsection{Updating the matrix \W}\label{sec:updateMatrixW}
The optimization problem in~(\ref{eq_als1}) is a least-squares problem with multiple right-hand sides. It can be rewritten as follows as a problem with a single right-hand side:
\begin{align}
	& \min_{\W} \normf{\Jcal^T_{(1)} - (\H\odot\V) \W^T}^2 \nonumber \\
	=& \min_{\W} \norm{\vec{\Jcal^T_{(1)}} - \vec{(\H\odot\V) \W^T}}^2 \nonumber \\
	=& \min_{\W} \norm{\vec{\Jcal^T_{(1)}} - \bigl(\I_n \otimes (\H\odot\V)\bigr)\vec{\W^T}}^2 \label{eq_lastMinEq}
\end{align}	
Here, $\I_n \in\ER^{n\times n}$ is the identity matrix and we denote the Kronecker product as $\otimes$. If we denote $\B_1 = \I_n \otimes (\H\odot\V)$, so
$$
\B_1 = 
\underbrace{\begin{bmatrix}
	\H\odot \V & \\
	&\ddots\\
	&&\H\odot \V 
\end{bmatrix}}_{\textrm{$n$ repetitions}}
$$
then the minimization problem in~(\ref{eq_lastMinEq}) can be rewritten as
\begin{equation}
\min_{\W} \norm{\vec{\Jcal_{(1)}^T} - \B_1 \vec{\W^T}}^2.
\label{eq:1rightHandSide}\end{equation}

On this level, the permuted weight $\weight_{(1)}$ is included as an extra factor in the cost function, namely
\begin{align}
&\min_{\W}  \norm{\vec{\Jcal_{(1)}^T} - \B_1 \vec{\W^T}}^2_{\weightone}  \label{eq_weightedapprox}\\
=& \min_{\W}  \bigl(\vec{\Jcal_{(1)}^T} - \B_1 \vec{\W^T} \bigr)^T \weightone \bigl(\vec{\Jcal_{(1)}^T} - \B_1 \vec{\W^T} \bigr). \nonumber
\end{align}
The weighted least squares solution is then given by
\begin{equation}
\vec{\W^T} = 
\left(\B_1^T \weightone \B_1\right)^{\!-1} \B_1^T \weightone\vec{\Jcal_{(1)}^T}.
\label{eqSolutionALS}\end{equation}

Here, care should be taken as the weight matrix \weight\ should be permuted to $\weightone$ using the same permutation that permutes $\vec{\Jcal}$ to
$\vec{\Jcal_{(1)}^T}$. If we denote the permutation matrix as $\P_1$, then we have
$$
\P_1 \vec{\Jcal} = \vec{\Jcal_{(1)}^T}
$$
and
$$
\weightone = \P_1 \weight\P_1^T.
$$
Using \matlab, the permutation matrix $\P_1$ is found by the following instructions:
\begin{verbatim}
   T = reshape(1:m*n*N, [n,m,N]); I = eye(m*n*N);
   P_1 = I(vec(reshape(T, [n, m*N])'), :); 
\end{verbatim}
Finally, reshaping the solution vector of \eqref{eqSolutionALS} to the dimensions of $\W^T$ gives after transposition the updated matrix $\W$.

\subsubsection{Updating the matrices \V\ and \H}
\begin{figure}[ht]
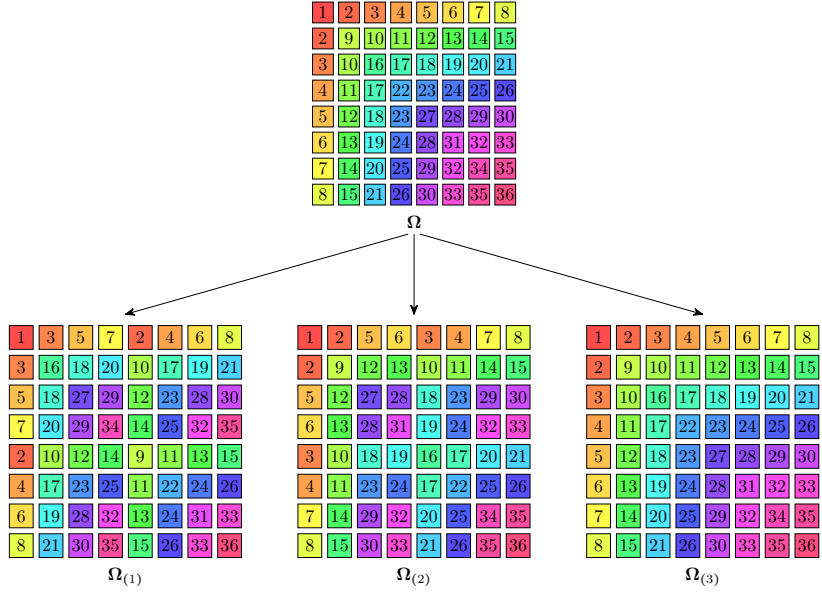

\pgfplotsset{compat=newest}
\usetikzlibrary{plotmarks}
\usepgfplotslibrary{patchplots}
\definecolor{mycolor1}{rgb}{0.90200,1.00000,0.30000}%
\definecolor{mycolor2}{rgb}{1.00000,0.98400,0.30000}%
\definecolor{mycolor3}{rgb}{1.00000,0.87000,0.30000}%
\definecolor{mycolor4}{rgb}{1.00000,0.75600,0.30000}%
\definecolor{mycolor5}{rgb}{1.00000,0.64200,0.30000}%
\definecolor{mycolor6}{rgb}{1.00000,0.52800,0.30000}%
\definecolor{mycolor7}{rgb}{1.00000,0.41400,0.30000}%
\definecolor{mycolor8}{rgb}{0.30000,1.00000,0.49600}%
\definecolor{mycolor9}{rgb}{0.30000,1.00000,0.38200}%
\definecolor{mycolor10}{rgb}{0.33200,1.00000,0.30000}%
\definecolor{mycolor11}{rgb}{0.30000,0.82000,1.00000}%
\definecolor{mycolor12}{rgb}{0.30000,0.93400,1.00000}%
\definecolor{mycolor13}{rgb}{0.30000,1.00000,0.95200}%
\definecolor{mycolor14}{rgb}{0.35000,0.30000,1.00000}%
\definecolor{mycolor15}{rgb}{0.30000,0.36400,1.00000}%
\definecolor{mycolor16}{rgb}{0.80600,0.30000,1.00000}%
\definecolor{mycolor17}{rgb}{1.00000,0.30000,0.85200}%
\definecolor{mycolor18}{rgb}{1.00000,0.30000,0.96600}%
\definecolor{mycolor19}{rgb}{0.92000,0.30000,1.00000}%
\definecolor{mycolor20}{rgb}{1.00000,0.30000,0.62400}%
\definecolor{mycolor21}{rgb}{1.00000,0.30000,0.51000}%

\begin{center}
\resizebox{0.99\textwidth}{!}{
}%
\end{center}
	\caption{A graphical representation of how the permutations of $\weight$ work, in the case where $m=n=N=2$. These permutations are used due to the vectorizations of several tensors in Section~\ref{sec:updateMatrixW} and beyond. Under the hood, these are defined in the same fashion as \matlab's {\tt reshape} function. We note that the permutation matrix $\P_3$ is, in fact, the identity matrix.}
	\label{fig_coloredGrids}
\end{figure}

When updating the matrix \V, the permutations from $\vec\Jcal$ to $\vec{\Jcal_{(2)}^T}$ and $\vec{\Jcal_{(3)}^T}$ should also be updated. This is done similarly as in Section~\ref{sec:updateMatrixW}. \fig{fig_coloredGrids} shows a graphical representation of the permuted matrices $\weightone$, $\weighttwo$ and $\weightthree$.

\subsection{Using dense weights}\label{sec:f} As mentioned earlier, the dense covariance matrix~$\mathbf{\Sigma}_\Jcal^\textrm{d}$ is rank-deficient and does not have an inverse, which could be used as weight matrix. Although the pseudo-inverse of~$\mathbf{\Sigma}_\Jcal$ could be computed, it does not yield enough equations in order to solve the problem~(\ref{eq:costfunction}) uniquely: the system~(\ref{eq:costfunction}) is in fact undetermined. That is why the following technique will be used in order to incorporate the weight matrix into the CPD decomposition. It consists of two parts: Section~\ref{sec:fWeight_part1} describes the first set of equations, and Section~\ref{sec:fWeight_part2} contains the second set of equations.

\subsubsection{Using the singular values of $\mathbf{\Sigma}_\Jcal$}\label{sec:fWeight_part1}
Let $\rSigma = \rank(\mathbf{\Sigma}_\Jcal)$ denote the rank of the covariance matrix, then we can decompose the matrix $\mathbf{\Sigma}_\Jcal$ using the singular value decomposition (SVD), and obtain
\setlength{\arrayrulewidth}{0.7pt}
$$
\rule[-15pt]{0pt}{35pt}
\mathbf{\Sigma}_\Jcal = \Vs\, \Ds\, \Vs^T = 
\Bigl[
\begin{array}{c|c}
 		\smash{\overbrace{\Vs^{(1)}}^{mnN\times \rSigma}}
 &
 		\smash{\underbrace{\Vs^{(2)}}_{\raisebox{-2pt}{\makebox[20pt]{$\scriptstyle mnN\times(mnN-\rSigma)$}}}}
\end{array}
 \Bigr] \,
\left[\begin{array}{@{}c|c@{}}
\rule[-5pt]{0pt}{0pt}\smash{\overbrace{\Ds^{(1)}}^{\rSigma\times\rSigma}}
 & 0 \\
\hline
0 & 0
\end{array} \right] \,
\Bigl[\begin{array}{@{}c|c@{}}
		\smash{\Vs^{(1)}}
 &
		\smash{\Vs^{(2)} \Bigr]^T}
\end{array}
=
\Vs^{(1)}\,  \Ds^{(1)}\, (\Vs^{(1)})^T.
$$
Here, $\Vs$ is an orthogonal matrix and~$\Ds$ is a diagonal matrix containing the singular values of $\mathbf{\Sigma}_\Jcal$. We call $\Vs^{(1)}$ the submatrix of $\Vs$ containing the first $\rSigma = \rank(\mathbf{\Sigma}_\J)$ columns of $\Vs$, and $\Ds^{(1)}$ the submatrix of $\Ds$ containing the non-zero singular values.

As in Section~\ref{sec:updateMatrixW}, every factor during the Alternating Least Squares algorithm uses a permuted version of $\mathbf{\Sigma}_\Jcal$. For example, at the update of the factor $\W$, we define
\begin{equation}
\mathbf{\Sigma}_{(1)} = \P_1  \mathbf{\Sigma}_\Jcal \P_1^T,
\label{eq_permCare}\end{equation}
where $\P_1$ is the permutation matrix, as defined in Section~\ref{sec:permutation}. It follows that the SVD of $\mathbf{\Sigma}_{(1)}$ is given by 
$$
	\mathbf{\Sigma}_{(1)} = \V_{\mathbf{\Sigma}_{(1)}} \D_{\mathbf{\Sigma}_{(1)}} \V_{\mathbf{\Sigma}_{(1)}}^T = (\P_1 \Vs^{(1)}) \,\Ds^{(1)} \, \bigl((\Vs^{(1)})^T\P_1^T\bigr)
$$

Let $\Q_1$ denote the $\rSigma\times mnN$-matrix
$$
	\Q_1 = \sqrt{(\Ds^1)^{-1}} (\Vs^{(1)})^T \P_1^T,
$$
such that $\Q_1^T \Q_1 = \weightone$. Then the solution of the minimization problem~(\ref{eq_weightedapprox}) is given by
\begin{equation}
(\Q_1  \B_1)^{\!\dagger} (\Q_1  \vec{\Jcal_{(1)}}^T),
\label{eq_lowRankALS}\end{equation}
where $\B_1$ is the block-matrix of Khatri-Rao products and $\X^\dagger$ denotes the pseudo-inverse of \X. In Appendix~\ref{sec:app1}, a derivation for solution~(\ref{eq_lowRankALS}) is shown. Similar expressions are found for the other factors being updated. The following table shows how the matrix~$\B_i$ changes, depending on the updated factor:

\medskip
\centerline{%
\begin{tabular}{p{0.5cm} p{1cm} c  p{2cm} p{2cm}}
$i$ & Updating factor & Matrix $\B_i$ & Dimensions of $\B_i$ & Dimensions of $\Q_i \B_i$\\
\hline
\rule{0pt}{7ex}1 & \W & 
$\begin{bmatrix}
\H\odot \V & \\
&\ddots\\
&&\H\odot \V 
\end{bmatrix}$
& $mnN\times rn$
& $\rSigma \times rn$
\\[2em]\hline
\rule{0pt}{7ex}2 & \V &
$\begin{bmatrix}
\H\odot \W & \\
&\ddots\\
&&\H\odot \W 
\end{bmatrix}$
& $mnN \times rm$
& $\rSigma \times rm$
\\[2em] \hline
\rule{0pt}{7ex}3 & \H &
$\begin{bmatrix}
\V\odot \W & \\
&\ddots\\
&&\V\odot \W 
\end{bmatrix}$
& $mnN\times rN$
& $\rSigma \times rN$
\end{tabular}}

\medskip
\noindent As we can see from this table, the dimensions of the systems being solved in (\ref{eq_lowRankALS}) depend on the factor being updated. We note that if the number~$N$ of operating points is too large, then the system with coefficient matrix $\Q_3\B_3$ is underdetermined, whereas the two other systems are not, as long as $\rSigma > \max\{rn,rm\}$. To remedy this, the~$\Vs^{(2)}$-part of the $\Vs$-matrix will be used to add extra constraints to the current systems of equations. These extra conditions will also be used for updating the matrices~\W\ and~\V.

\subsubsection{Adding more equations to the existing set}\label{sec:fWeight_part2}
If the covariance matrix~$\mathbf{\Sigma}_\Jcal$ of the noise $\v$ is not of full rank, there exist linear relations $\L$ between the noise disturbances $\v$, such that $\L \v = 0$. In Appendix~\ref{sec:app2}, it is shown how these can be retrieved from $\mathbf{\Sigma}_\Jcal$.

 In order to find the extra equations for the update of the current factor, we reconsider the minimization problem~(\ref{eq_weightedapprox}) for the updating factor $\W$:
$$
\min_{\textrm{vec}({\W^T})} \norm{\vec{\Jcal_{(1)}^T} - \B_1  \vec{\W^T}}^2_{\weightone}
$$
Rewriting this expression as $\vec{\Jcal_{(1)}^T} = \B_1  \vec{\W^T} + \v$, where $\v$ is correlated noise with $\mathbf{\Sigma}_{(1)}$ as in~\eqref{eq_permCare}, Appendix~\ref{sec:app2} can be used in order to add extra equations to the existing set in~(\ref{eq_lowRankALS}). We remark that the noise $\v$ is correlated, due to the linear relations between the covariance matrices $\mathbf{\Sigma}_f$ and $\mathbf{\Sigma}_\Jcal$. 

Using the notations in Equation~(\ref{eq_permCare}), the orthogonal factor of the singular value decomposition of $\mathbf{\Sigma}_{(1)}$ is given by
$$
\P_1 \Vs = \Bigl[ 
\rule[-20pt]{0pt}{30pt}
\begin{array}{@{}c|c@{}}
		\smash{\overbrace{\P_1\Vs^{(1)}}^{\makebox[0pt]{$\scriptstyle mnN\times \rSigma$}}}
 &
 		\smash{\underbrace{\P_1\Vs^{(2)}}_{\makebox[0pt]{$\scriptstyle mnN\times(mnN-\rSigma$)}}}
\end{array}
\Bigr].
$$
In Appendix~\ref{sec:app2}, it is proven that the extra conditions are then given by
\begin{equation}
\big( (\Vs^{(2)})^T \P_1^T \B_1 \big)^{\!\dagger}\, \bigl((\Vs^{(2)})^T \P_1^T \vec{\Jcal_{(1)}^T} \bigr).
\label{eq_resultEquation2}\end{equation}
Similar expression can be found for the other two updating factors.

In general, we conclude that the following system of equations must be solved for every updated factor $1\le i\le 3$:
\begin{equation}
\begin{bmatrix}
\Q_i \B_i \\
(\U_\Jcal^2)^T \P_i^T \B_i 
\end{bmatrix} 
^{\!\dagger}
\begin{bmatrix}
\Q_i \vec{\Jcal_{(i)}}^T \\
(\U_\Jcal^2)^T \P_i^T \vec{\Jcal_{(i)}^T}
\end{bmatrix}.
\label{eq_totalEq}\end{equation}

\subsection{Stopping criteria for the \CPD\ algorithm}\label{sec_stopcriteria}
By iterating the presented algorithm, a sequence of updated factors is obtained
$$
\W^{(1)}, \V^{(1)}, \H^{(1)}, \W^{(2)}, \V^{(2)}, \H^{(2)}, \W^{(3)}, \V^{(3)}, \H^{(3)}, \ldots
$$
We have the following two stopping criteria for this iteration algorithm:
\begin{enumerate}
\item When the relative step size between two iterations is below a given tolerance, then the algorithm is stopped. The relative step size at iteration step $j\ge 2$ is given as
$$
\frac{\sqrt{ \normf{\W^{(j)} - \W^{(j-1)}}^2 + \normf{\V^{(j)} - \V^{(j-1)}}^2 + \normf{\H^{(j)} - \H^{(j-1)}}^2}}
{\sqrt{ \normf{\W^{(j)}}^2 + \normf{\V^{(j)}}^2 + \normf{\H^{(j)}}^2 }}\textrm{, or}
$$
\item when an upper bound on the number of iterations is reached, then the algorithm is stopped. This takes care of possible divergence.
\end{enumerate}

\subsection{Summary of the proposed algorithm}
In this section, we summarize the proposed weighted \CPD\ algorithm and incorporate it into the larger decoupling problem of the noisy multivariate polynomial function~\f. We assume the covariance matrix~$\mathbf{\Sigma}_\f$ of the coefficients of~\f\ to be known and the same notations as in Section~\ref{sec:intro} will be used.

\begin{algor}{Decomposition of the noisy multivariate polynomial \f, given $\mathbf{\Sigma}_\f$.}
\begin{enumerate}
	\item Evaluate the Jacobian matrix $\J(\u)$ of \f\ in $N$ randomly chosen sampling points $\u^{(1)},\ldots,\u^{(N)}$.
	\item Stack the Jacobians into a three-way tensor $\Jcal$.
	\item Transform the covariance matrix $\mathbf{\Sigma}_\f$ of \f\ into the matrix~$\mathbf{\Sigma}_\Jcal$. Here, three choices are possible: the element-wise as discussed in~(\ref{eq_elementwiseCovar}), slice-wise as discussed in~(\ref{eq_slicewise}) or dense covariance matrix discussed in~(\ref{eq_densecovar}).
	\item Compute the Weighted Canonical Polyadic Decomposition (WCPD) of $\Jcal \approx \sum_{i=1}^r \w_i\circ\v_i\circ\h_i$. In the element- and slice-wise case, use~(\ref{eqSolutionALS}), in the dense case, use~(\ref{eq_totalEq}). Iterate until one of the stopping criteria in Section~\ref{sec_stopcriteria} is satisfied.
	\item Starting from the vectors $\h_1,\ldots,\h_r$, reconstruct the internal univariate functions $g_1(x_1),\ldots,g_r(x_r)$ with an integration step, as in Algorithm~1.
\end{enumerate}
\end{algor}

\section{Numerical experiments}\label{sec:results}
In this section, results of the methods suggested in this paper will be discussed. 

\subsection{Correlations of the errors of the CPD}
In this section, we demonstrate that the weighted CPD decomposition works as expected. For this, we set $m = n = N = 2$ and start with a given $8\times8$ positive definite weight matrix
\begin{equation}
\weight = \begin{bmatrix}                                            
   0.74 & 0 & 0 & 0 & 0 & 0 & 0 & 0 \\
   0 & 1.67 & 0 & 0 & \color{blue}\underline{0.87} & 0 & 0 & 0 \\
   0 & 0 & 0.96 & 0 & 0 & 0 & 0 & \color{red}\underline{\underline{0}} \\
   0 & 0 & 0 & 0.63 & 0 & 0 & 0 & 0 \\
   0 & \color{blue}\underline{0.87} & 0 & 0 & 1 & 0 & 0 & 0 \\
   0 & 0 & 0 & 0 & 0 & 0.11 & 0 & 0 \\
   0 & 0 & 0 & 0 & 0 & 0 & 0.77 & 0 \\
   0 & 0 & \color{red}\underline{\underline{0}} & 0 & 0 & 0 & 0 & 0.31 \\
\end{bmatrix}.
\label{eq_wDefinition}\end{equation}
This matrix is chosen to be almost diagonal, and has two extra nonzero covariance elements. The tensor \Tcal\ to be decomposed has eight elements indexed as follows: 
\begin{equation}
\Tcal_1 = \begin{bmatrix}
t_1 & t_3 \\
t_2 & t_4
\end{bmatrix} \quadtext{and} \Tcal_2 = \begin{bmatrix}
t_5 & t_7 \\
t_6 & t_8
\end{bmatrix},
\label{eq_tensorNumbering}\end{equation}
where $\Tcal_1$ and $\Tcal_2$ denote the first and second frontal slices of \Tcal\ respectively. Elements $t_2$ and $t_5$ are correlated with cross-covariance 0.87, which is shown in {\color{blue}\underline{blue}} in Equation~\ref{eq_wDefinition}, while elements $t_3$ and~$t_8$ are uncorrelated, and are shown in {\color{red}\underline{\underline{red}}}.

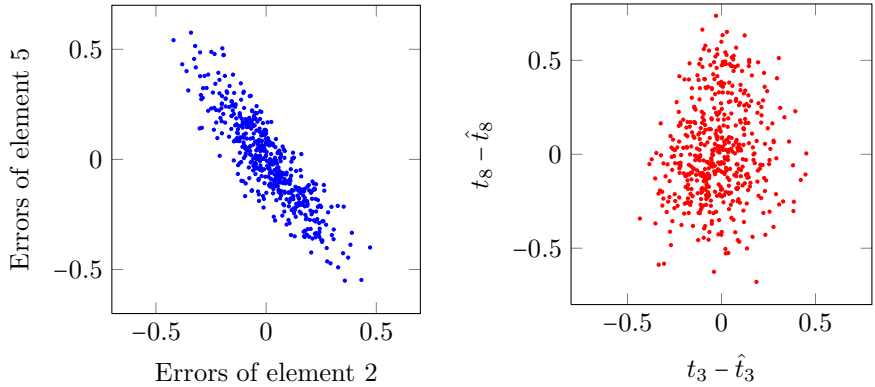
\begin{figure}[htpb]
\centerline{
\begin{tikzpicture}
    \begin{axis}[
        xlabel={Errors of element 2},
        ylabel={Errors of element 5},
        xmin = -0.7,
        xmax = 0.7,
        ymin = -0.7,
        ymax = 0.7,
        xtick = {-0.5, 0, 0.5},
        ytick = {-0.5, 0, 0.5},
        width = 5 cm,
        height = 5 cm,
        scale = 1.2,
    ]
\addplot[only marks, blue, mark options={scale  =0.3}] coordinates {
(0.1871708300, -0.2104278800)
(-0.1133672700, 0.2712248100)
(-0.0231752060, -0.0399206660)
(-0.1067146900, 0.1445062200)
(-0.0217063950, 0.2156476700)
(-0.0468389660, -0.0184808920)
(0.2552159700, -0.2810592200)
(0.2213955500, -0.3374314200)
(0.0627220810, 0.0057839655)
(0.0683935920, -0.1435437900)
(0.2420615800, -0.3295179100)
(0.2965722400, -0.2152670000)
(0.1614415200, -0.1385655100)
(0.2012774000, -0.2324770800)
(0.0831064340, 0.0071358786)
(-0.3266787900, 0.4564125100)
(0.2271715800, -0.2363545800)
(-0.1316724800, 0.1217476700)
(0.1759335600, -0.2130020700)
(-0.0772270740, 0.1608152800)
(0.2495086900, -0.2157809400)
(-0.1148014600, 0.3455875000)
(-0.0983328070, 0.0561291500)
(0.0363384140, -0.0949159250)
(-0.0484862250, -0.0662008690)
(-0.0637272160, 0.2224114700)
(0.2155434400, -0.1093352300)
(-0.1400982000, 0.1546646700)
(-0.2273365900, 0.2246960900)
(0.1157255600, -0.2299300900)
(0.1771245400, -0.3236883000)
(-0.1576313400, 0.1975735000)
(0.0402612900, -0.0020207519)
(0.0634246630, -0.1531805700)
(-0.0757879310, 0.0565013440)
(0.2037242600, -0.1515421900)
(0.2451429000, -0.2101728200)
(-0.0495601260, 0.2118069500)
(0.0561844050, -0.0364823380)
(-0.2490754100, 0.2312012900)
(0.0000850798, 0.0713977540)
(-0.0040209267, 0.0109134850)
(-0.1986674600, 0.5045932500)
(0.1592161000, -0.0771462170)
(0.0454545540, -0.1164159300)
(-0.0639009600, -0.1403339900)
(-0.0272403540, 0.0535266370)
(0.2536992700, -0.3744326800)
(0.1013246900, -0.1508735800)
(-0.2795700700, 0.3790564500)
(0.0328265540, -0.1342093000)
(-0.0073992044, -0.1078826700)
(-0.2016270800, 0.2507246600)
(0.1525420400, -0.2374000900)
(0.2721864300, -0.2967554200)
(0.0659681170, -0.0082604211)
(-0.1461685300, 0.0821637050)
(-0.0832792160, 0.0476313730)
(0.2173009400, -0.2764130200)
(-0.2630617500, 0.2208296100)
(0.0103323580, 0.0862457890)
(0.0602627570, 0.0403458300)
(-0.1206592100, 0.0830886240)
(-0.0051981820, 0.1227693200)
(0.3723022900, -0.4458144000)
(0.0841079700, -0.1517019300)
(-0.1110766300, 0.1850617600)
(0.2349106200, -0.2154576900)
(0.2320763500, -0.2708961600)
(0.3068321600, -0.3182499300)
(-0.1352841600, 0.1860567500)
(-0.1839848900, 0.2764151700)
(-0.2116970600, 0.3000268500)
(-0.0745328960, 0.2870863600)
(-0.1281449700, -0.0659570410)
(0.0065328320, 0.1683598400)
(-0.1154362000, 0.0991573710)
(-0.1381683800, 0.2500581800)
(0.4716695100, -0.3998407500)
(-0.0827412180, 0.0434212150)
(-0.0117003230, -0.2184951100)
(-0.0935095660, 0.1827706500)
(0.0327568780, 0.0632546220)
(0.0131913520, -0.0109571930)
(-0.0322894460, -0.1883901900)
(0.0183421350, -0.1092157600)
(0.0382200920, -0.1226346600)
(-0.0943833620, 0.0580735120)
(0.0405920730, 0.0076372432)
(-0.0146062520, 0.1649167600)
(-0.0330585820, -0.1590827300)
(0.2081896000, -0.1787944000)
(0.0163223120, -0.0241611270)
(-0.1325187600, 0.0966584490)
(0.3826652500, -0.3875775800)
(-0.0368085620, -0.1452976900)
(-0.0559975200, 0.1177929300)
(0.1717904900, -0.1454572600)
(-0.1132723800, -0.0118712140)
(-0.0037846995, -0.0785253180)
(-0.2167813700, 0.3936508200)
(0.2116941100, -0.1408195900)
(0.0315819380, -0.1541265200)
(0.1723184000, -0.1414504400)
(0.1125913800, -0.1492326900)
(0.2502773000, -0.3652611000)
(-0.0710438580, 0.0943728730)
(0.1833884300, -0.2667064900)
(-0.2371166600, 0.1349115900)
(0.1545418600, -0.2070666400)
(0.0560650810, -0.1377618700)
(-0.1182580200, 0.0991539780)
(0.0609082310, -0.0533962140)
(0.2420359500, -0.2104437100)
(-0.2072005000, 0.1882745500)
(-0.0279394630, 0.0013616987)
(0.3520687900, -0.2144223600)
(-0.1192584400, 0.0302374950)
(0.0500347520, -0.0376022300)
(0.0009899749, -0.1515158200)
(-0.0478640730, 0.1847680100)
(0.1167813400, -0.1393354100)
(-0.0190611360, -0.0509279170)
(-0.1337340900, 0.1599525500)
(-0.0495494040, 0.0806680200)
(0.1029924300, -0.2129269300)
(-0.0360539420, 0.0662491610)
(-0.0952340730, -0.1668001200)
(-0.1462877500, 0.2060721800)
(-0.0136666190, -0.1662032700)
(-0.1001505100, 0.0989802380)
(0.3061407200, -0.3623792400)
(-0.0581129430, 0.0666230190)
(0.2214085700, -0.1362917400)
(-0.0857478540, 0.1543580000)
(-0.2264595500, 0.1675198600)
(0.2249666800, -0.3535418700)
(-0.2959887700, 0.2923184900)
(0.2983092100, -0.3439299000)
(-0.0921878480, 0.1175873200)
(-0.1241920900, 0.1503180600)
(-0.1402729800, 0.2391053500)
(0.1759220500, -0.1329629500)
(0.1205541800, 0.0052410842)
(-0.0944371850, 0.0495334740)
(0.0776780340, -0.2028653700)
(-0.0271238450, 0.0830849810)
(-0.1943013500, 0.2064036600)
(-0.0379226880, 0.0681281400)
(0.1856979800, -0.1913810200)
(0.0184862870, 0.0772297820)
(0.0494998140, -0.0188098300)
(0.1377951600, -0.2128942700)
(0.2271345200, -0.3348814300)
(0.0235928930, 0.1701054300)
(-0.0776197610, 0.1797510700)
(-0.2915249600, 0.1436085800)
(-0.1450323000, 0.3021111800)
(0.0265818370, -0.0205516230)
(-0.0786908570, -0.0365792650)
(-0.0215617080, 0.0017057955)
(0.1701338400, -0.3281770800)
(-0.1335003300, 0.0836245830)
(-0.0998136290, 0.2337137000)
(-0.2996103800, 0.4860475500)
(-0.0713080920, 0.0601763300)
(0.1358907500, -0.0525333070)
(-0.0013826570, 0.0930436660)
(0.2430137700, -0.3143056000)
(0.0187758120, -0.1197981800)
(-0.1676843900, 0.2374948700)
(0.0687831650, -0.1040871900)
(0.1472972600, -0.1736396600)
(-0.0272171600, -0.0464810720)
(-0.0128569240, 0.1173168600)
(-0.0696222300, 0.2889359800)
(0.2770881500, -0.4632004400)
(0.0781023550, -0.0665654740)
(-0.2001077400, 0.2875243300)
(0.2184337000, -0.1324195200)
(-0.0321565770, 0.0889107220)
(-0.0745574200, 0.0245641290)
(-0.0578046240, -0.0294829870)
(-0.1196519500, 0.2097758800)
(0.0610454420, -0.0592005580)
(0.1755879200, -0.1746073700)
(-0.0070820081, 0.1344041300)
(-0.1871471900, 0.3793181200)
(-0.0635334460, 0.2186267000)
(0.2376019000, -0.3401382900)
(-0.0094244761, -0.0447933970)
(0.0845847660, 0.0381175340)
(0.0806521660, 0.0451576120)
(0.1451238300, -0.1311008700)
(0.1126165500, -0.3434360900)
(0.1262358900, 0.0078639420)
(-0.1765048400, 0.3021572000)
(-0.0577463580, 0.1520127100)
(-0.1145070100, 0.1536138700)
(0.1199699700, 0.0031212951)
(0.0405833880, -0.0966890370)
(-0.1405139700, 0.0868305090)
(-0.1793513800, 0.3090134400)
(-0.1042859500, -0.0985489670)
(-0.0009558862, 0.1301052200)
(0.0370109820, -0.1071946200)
(-0.0180909850, 0.0843849480)
(0.0419039420, 0.0558297650)
(0.0928926330, -0.1814709600)
(0.2921386500, -0.4718841700)
(0.1955750900, -0.3511895200)
(-0.0077773751, 0.0485246050)
(-0.0482340160, -0.0419860690)
(-0.1119710100, -0.0971283500)
(-0.2483698100, 0.4882129600)
(-0.0278176260, 0.0172073170)
(0.2105118900, -0.3594668600)
(0.0653216770, -0.0802607780)
(-0.1194378100, -0.0161743290)
(0.1344670000, -0.2605184000)
(0.0571067530, -0.1047442400)
(0.2749068600, -0.4023786200)
(-0.1000489400, 0.0641057130)
(0.0670549870, -0.1024576400)
(-0.2212999800, 0.0967565340)
(0.3570563000, -0.5504220400)
(0.0928916720, -0.2956284700)
(-0.0005304260, 0.0773398750)
(-0.2859915400, 0.3238124500)
(-0.0449204250, 0.0014250687)
(0.1882125600, -0.2019377500)
(0.2353035600, -0.1848396400)
(0.1262472300, -0.2445514500)
(0.1730230600, -0.1899655600)
(-0.1560301800, 0.1971312500)
(-0.0921806150, 0.0266820450)
(0.0016652467, 0.0331658910)
(-0.4202177000, 0.5415623500)
(0.0406778100, -0.2774790300)
(0.4323932600, -0.5473736400)
(0.0377068680, -0.0487372770)
(0.0610933410, -0.0931764330)
(-0.0439655190, 0.1313728400)
(-0.1545475100, 0.2447077900)
(-0.2338352900, 0.4786754400)
(-0.0084809497, 0.0386011360)
(0.0463847880, -0.1851826100)
(0.0275987720, -0.0912487810)
(0.1119377700, -0.0444489930)
(0.1684189300, -0.1640813800)
(-0.1277345300, 0.1368962700)
(0.0014884522, 0.0357747510)
(0.0457248060, -0.0053687613)
(-0.1924729800, 0.4743728700)
(-0.3611347200, 0.4009143700)
(0.0196106070, -0.0871096930)
(0.2125794300, -0.2380304300)
(-0.0420827130, 0.2098035500)
(-0.0945609880, 0.0158450960)
(-0.1022788200, 0.0770605770)
(0.1936921000, -0.1165437800)
(-0.2118142000, 0.3611146000)
(0.1560321900, -0.1523195300)
(-0.0763155390, 0.1886420000)
(0.0865280010, -0.0655633610)
(-0.0278220480, 0.1247277800)
(-0.1184561400, 0.1951284500)
(0.0749246200, 0.0586644300)
(0.1131601300, -0.1804959600)
(-0.1581610400, 0.1385694500)
(-0.1469354900, -0.0045763291)
(0.0894534860, -0.1089859200)
(-0.2029784200, 0.2320017500)
(0.0825472540, -0.0886143850)
(-0.1836522300, 0.2405429300)
(0.0289614880, 0.0676138360)
(0.0635870750, -0.1766680900)
(0.0934296840, -0.0600938750)
(-0.0343726200, -0.0939235130)
(0.0444707600, 0.0926456190)
(0.0363598110, -0.0188464530)
(0.0228783270, -0.0027474532)
(-0.0269311360, 0.2182768300)
(-0.1280629600, 0.1180786400)
(0.1807935800, -0.2714589900)
(-0.2629582300, 0.3852650100)
(0.0358878890, -0.1129202500)
(-0.1562489200, 0.2642295200)
(0.1233613300, -0.1666552300)
(0.2251860600, -0.1797003800)
(0.2551463700, -0.3683824000)
(0.2059243800, -0.2490035300)
(-0.0346706720, 0.0524510500)
(-0.0661058150, 0.1713025300)
(-0.3015754600, 0.1408390600)
(0.0131271950, -0.1171453000)
(0.1234315200, -0.1809534400)
(0.0869286050, -0.0957865680)
(-0.0104400970, -0.1793849800)
(0.2639856700, -0.1417849800)
(-0.0588100230, -0.0499852680)
(-0.0725146680, 0.0384847590)
(-0.1485861900, 0.2495443600)
(0.0084095079, 0.0963683420)
(-0.3193665000, 0.4178124900)
(0.1628215700, -0.1985483000)
(0.2041563400, -0.1477047000)
(-0.0185239420, 0.0651938240)
(0.0519252680, 0.0290576710)
(0.1136887100, -0.0737988230)
(-0.1199985900, 0.2017467000)
(-0.1934723700, 0.4735365000)
(-0.0384507770, 0.0491771750)
(-0.0003026797, 0.1178049700)
(0.0455673280, -0.0827189010)
(0.0042815225, -0.0075670264)
(-0.3409571700, 0.5758110600)
(0.1780671100, -0.2518925500)
(0.3269028300, -0.4901806800)
(-0.2940659900, 0.2732728300)
(0.1838733800, -0.3412749600)
(-0.0274189720, 0.1837481700)
(-0.2137458300, 0.2220316300)
(-0.0677062610, 0.1396807900)
(0.1306206200, -0.2289095800)
(-0.2302699800, 0.3435277100)
(0.1403868500, -0.2133235800)
(0.0087254574, -0.0395904710)
(-0.1271809300, -0.0259941350)
(0.1860071800, -0.1896674200)
(0.0382123810, -0.0684637390)
(0.1491633900, -0.1951967700)
(0.0734473300, -0.2551930000)
(-0.1009551100, 0.0059806475)
(0.0552420090, 0.0358650260)
(-0.1659339500, 0.3336340800)
(-0.1503571100, 0.3849911800)
(0.0529864350, -0.1220430100)
(0.3249846200, -0.2111096800)
(0.0177480450, -0.0063136227)
(-0.0323033250, -0.0438090670)
(0.1101620100, -0.0971821540)
(0.1430802300, -0.1334156100)
(0.0643322220, -0.1038759300)
(0.1749052400, -0.3523726800)
(-0.3806935500, 0.4319434500)
(0.3890301900, -0.3330746500)
(0.1107424000, -0.1773478200)
(-0.1074995400, 0.1554799200)
(-0.0687544850, 0.0201969610)
(0.1220067500, 0.0304965110)
(0.0139656440, 0.0634494640)
(0.1625352700, -0.1373530400)
(0.0505259570, -0.0570853720)
(-0.0220085980, 0.1084118900)
(0.0785282230, -0.0631624750)
(0.2414992700, -0.3517516300)
(0.1557729100, -0.0487398300)
(-0.0261123040, 0.0378914180)
(0.2599053200, -0.2806756100)
(-0.0768954420, 0.1476256700)
(-0.0659781450, -0.1426982000)
(-0.0018270909, 0.1326733100)
(0.0561968920, -0.0999261510)
(0.1147284600, -0.0101269110)
(-0.0220374050, 0.0408490530)
(0.1332668900, -0.1244154800)
(-0.0063957822, 0.1103606400)
(0.0347175190, -0.2602206900)
(0.0051235637, -0.0289688810)
(-0.2700935300, 0.3243523500)
(-0.0599456710, 0.0972416960)
(0.1873699800, -0.1237590800)
(-0.1015321400, 0.3076207700)
(0.1022735300, -0.1496007900)
(-0.1013583500, 0.1126418400)
(-0.0520319090, 0.1548047800)
(0.0300095680, 0.0716180760)
(0.1712491600, -0.1894713300)
(-0.2006911300, 0.1185347200)
(0.0168003050, -0.0623357330)
(-0.0258611180, -0.0275053380)
(-0.0115476180, 0.0581868940)
(-0.0297722320, 0.1399667100)
(0.0916410610, -0.0686623610)
(0.0666375760, 0.0252970910)
(-0.0585548470, 0.1827685000)
(0.1056194800, -0.1937033800)
(0.2101669300, -0.1000946600)
(0.1237060600, -0.1412735600)
(0.0165546770, -0.0473283270)
(0.0555217680, -0.0831046000)
(-0.0317588570, 0.1287288600)
(0.1289753800, -0.2274377500)
(0.1986109300, -0.2746681100)
(0.1243273900, -0.2225688100)
(0.2267364200, -0.1125245500)
(0.0698078960, -0.0734781130)
(-0.0006618566, 0.0176347030)
(0.1000899600, -0.0984524720)
(-0.0498223760, -0.0166216520)
(-0.0435125820, 0.0816117910)
(-0.0152042710, 0.1108716600)
(-0.0874822180, 0.1689987400)
(-0.0500770430, -0.0181899430)
(0.0696985610, -0.0599295490)
(-0.1681879900, 0.1032735900)
(-0.0007361817, 0.1186337400)
(-0.0580530260, 0.0050174757)
(0.0816889370, -0.2084002500)
(-0.0699241650, -0.0138154990)
(0.1494874000, -0.1781350000)
(0.1679364900, -0.2065696100)
(-0.0879454190, 0.0947379160)
(-0.2020770300, 0.1736903700)
(-0.0353474270, 0.0137380130)
(0.1014184800, -0.2751904300)
(-0.0762350820, -0.0852094810)
(0.1062445700, -0.2897746700)
(0.0191239430, -0.1234854300)
(0.0246198420, -0.0597362010)
(-0.2994423900, 0.3778242200)
(-0.1860752800, 0.1671574800)
(0.1295303400, -0.3491635100)
(0.3492228400, -0.4261168600)
(-0.1832257900, 0.2754535100)
(-0.0015759289, 0.0684561500)
(0.0287645270, -0.0230083100)
(0.1373313700, -0.1350970600)
(0.1598746000, 0.0242323670)
(-0.1567651300, 0.2631463400)
(0.2618185400, -0.1159815900)
(0.1886186000, -0.3167430200)
(-0.0927547780, -0.1091489800)
(0.0677205470, -0.0679356030)
(-0.0423896810, -0.0367818150)
(0.0306188770, -0.0165803500)
(0.1125845700, -0.1526079600)
(0.0683071240, -0.1630524100)
(-0.1515703000, 0.2220484900)
(-0.1070487400, 0.0999680740)
(-0.0477571270, 0.0191123380)
(0.0825888910, -0.2369332800)
(-0.0646877730, -0.0090420660)
(-0.1057393600, 0.1689390600)
(0.0540709170, -0.0699229710)
(0.2047162800, -0.4225935500)
(-0.2416289900, 0.3036058100)
(0.1469944400, -0.1423258500)
(0.1088853200, -0.0884809350)
(0.2020237700, -0.3155710200)
(0.0834827820, -0.1460955300)
(0.1952842800, -0.1406092200)
(0.2810693100, -0.2404347500)
(0.1490471000, -0.1303789100)
(-0.0153378110, 0.1241311300)
(-0.0548167560, 0.1978493800)
(-0.2004909900, 0.2541723600)
(-0.0207046460, 0.0427111080)
(-0.0085383603, 0.0024496490)
(-0.1075422000, 0.0570399800)
(-0.2107965500, 0.0162894650)
(0.0668989610, 0.0169697050)
(-0.0339410310, -0.0890016190)
(-0.0444790270, 0.0409278600)
(-0.1303152100, 0.0897572650)
(-0.0883893770, 0.1465922600)
(-0.2877826800, 0.2795185500)
(-0.0537078370, 0.0084104560)
(0.1721375700, -0.3485713600)
(-0.0941832750, 0.1625028800)
(-0.3530944600, 0.3131799700)
(0.1590195900, -0.2454339300)
(-0.0208942020, -0.0755687190)
(-0.0885346730, 0.0419067390)
(-0.0590257010, 0.1569290700)
(0.1065168100, -0.1112645500)
(0.2245321800, -0.2817995700)
(0.0779631700, -0.1073644900)
(0.1426690300, -0.1641370100)
(-0.0491668020, -0.0223965030)
(0.0570353030, -0.1475359900)
(-0.0277856700, -0.0081252017)
(0.0001362934, 0.0713080910)
(0.0142028250, -0.0384212070)
(-0.3219444600, 0.5150612200)
(-0.2355013100, 0.2344315000)
(-0.0318479230, -0.0018339998)
(0.0106826130, -0.0784857330)
(0.0683870910, -0.0780666150)
(-0.0113574180, 0.0579880830)
(-0.0397117770, 0.1716336300)
(0.0227978460, 0.0243051730)
(0.1365505200, -0.2019501700)
(0.0177236500, 0.0011840714)
(-0.0151817340, 0.0036385720)
(-0.0362860480, 0.0788873570)
(-0.0434727520, 0.2422046000)
(-0.1032586700, 0.3021724900)
(0.2406321900, -0.2889192000)
};
\end{axis}
\end{tikzpicture}
\kern1em
\begin{tikzpicture}
    \begin{axis}[
        ylabel shift = -5 pt,
        xlabel={$t_3 - \hat{t}_3$},
        ylabel={$t_8 - \hat{t}_8$},
        xmin = -0.8,
        xmax = 0.8,
        ymin = -0.8,
        ymax = 0.8,
        xtick = {-0.5, 0, 0.5},
        ytick = {-0.5, 0, 0.5},
        width = 5 cm,
        height = 5 cm,
        scale only axis,
        scale = 0.8
    ]
\addplot[only marks, red, mark options={scale = 0.3}]  coordinates {
(-0.2129236200, -0.3763295300)
(0.0398952370, -0.2921205200)
(-0.0517627290, 0.3520956300)
(-0.1792527500, -0.1061516800)
(0.3803823100, -0.0588616480)
(-0.0361776020, 0.0288926560)
(-0.0533000550, -0.0312132370)
(-0.0227095920, 0.2096240600)
(0.0496092920, -0.2684978500)
(-0.0657535040, -0.1269985600)
(0.0584424470, 0.3643568200)
(0.2779286900, 0.0640301740)
(0.0022275579, -0.0946614140)
(-0.0384412000, -0.0908164650)
(0.1429474500, 0.1665925000)
(-0.1160892600, -0.0329105520)
(0.0890437210, -0.1552213700)
(0.0008852745, 0.4737887800)
(-0.0078768333, -0.0527500590)
(0.1263824700, 0.3390220000)
(0.0402356740, -0.2106098500)
(0.0639785180, -0.3490598900)
(-0.1908753000, -0.1740878100)
(0.0070169647, 0.3290009400)
(-0.0336608890, 0.3203785700)
(-0.2206454100, -0.0679404700)
(0.1649961000, -0.4393983100)
(-0.2388381300, 0.2774433900)
(-0.0458545110, 0.0163443750)
(-0.0259012440, 0.0914274550)
(-0.2238428800, -0.0748296820)
(-0.1697175500, 0.4839065800)
(-0.0565034280, 0.1336354700)
(-0.1967136200, 0.3984401200)
(-0.1157851700, 0.0812564740)
(-0.1666006000, -0.0007028927)
(0.1754674300, 0.0043540881)
(0.2183682900, -0.0573322650)
(0.0122425580, -0.0486550630)
(-0.1832868900, 0.0352257610)
(0.2042519500, -0.0760548050)
(0.1452390700, 0.4695835600)
(0.3279106800, -0.2510856800)
(0.1509153200, -0.1715104100)
(-0.0641205690, -0.4139020500)
(-0.2226831900, -0.0988111010)
(-0.0205877380, 0.0935142200)
(0.1233780100, 0.6324589000)
(-0.0181136220, -0.0795951980)
(-0.0260663980, 0.4599109300)
(-0.0319804930, 0.5205299100)
(-0.0266927360, 0.1264230900)
(-0.1984140500, -0.3154315900)
(0.0074700057, 0.3592299400)
(0.0104123270, -0.2555835900)
(0.0638509410, 0.0484398100)
(-0.0565245670, -0.2465785800)
(-0.1143164300, -0.0124511370)
(0.0414915460, -0.1373460400)
(-0.1059887500, -0.1857691100)
(0.2533031500, -0.0523952380)
(0.0356959000, -0.2040632900)
(-0.0518156830, 0.1861703700)
(-0.2766836300, -0.2279694600)
(-0.0074123485, -0.0873232980)
(0.0148206290, 0.3766387700)
(-0.1033997000, -0.2003209900)
(-0.0943333400, 0.1139394400)
(0.1191998800, 0.1328759700)
(0.0328118480, -0.0130146150)
(-0.0148779910, 0.3775383300)
(0.0343336920, -0.0303587030)
(-0.1168776900, 0.0242505650)
(0.1355110700, -0.5013665400)
(-0.2576799100, 0.0276275490)
(0.1123372600, -0.0662396630)
(-0.3832719600, -0.0523765740)
(-0.1268119800, -0.2074856000)
(-0.2460399100, -0.2486156500)
(-0.2177200400, 0.1577432000)
(-0.2183076500, 0.1692616700)
(0.0289291430, 0.6513347600)
(0.1700969800, -0.0301239400)
(-0.1247386200, -0.0249793950)
(-0.1016454800, 0.6633782800)
(-0.1569533700, -0.0976521970)
(-0.1784782100, -0.2350044000)
(-0.1596157300, -0.1339429200)
(0.2948990200, 0.4045148100)
(0.2264471800, 0.4327301800)
(-0.0352042540, 0.0546066980)
(0.0508104540, -0.2042876000)
(0.0217865250, -0.2007204300)
(-0.0015136394, 0.3856945400)
(0.1443665900, -0.1218992600)
(-0.2710992500, 0.0968939970)
(0.1428162600, 0.3140010700)
(0.0013543835, -0.4132068900)
(-0.0906917330, 0.4009735700)
(-0.0871523310, 0.1559514400)
(-0.0477417510, -0.0660713690)
(-0.0398203850, -0.6254806200)
(-0.2658413500, 0.0359578290)
(-0.2739225900, -0.4120002300)
(-0.1745004600, -0.3124022200)
(0.0201092390, 0.4939685100)
(0.1177466300, 0.1375146600)
(-0.0809077540, -0.2376090400)
(-0.0943598040, 0.5598832400)
(0.0136240170, -0.2177512600)
(0.0683737980, -0.1104667700)
(0.0235721650, 0.4805168800)
(0.0802561910, -0.0536217160)
(-0.0458661850, -0.0490990790)
(-0.0822203920, -0.2177908300)
(-0.1807565200, -0.1393479000)
(0.2902305400, -0.2124097100)
(-0.2759110700, -0.1429591100)
(0.1170788600, 0.0501629010)
(-0.1381553100, 0.2378840100)
(0.2598377800, 0.0287258010)
(0.0759001550, 0.0641493930)
(-0.0237764430, 0.5594662900)
(0.0130237030, 0.1849434800)
(0.0627396030, 0.5049866400)
(-0.0501260180, 0.2946959100)
(0.1626556900, 0.1100601200)
(-0.3220758700, -0.0056186973)
(0.0161492010, -0.0155910420)
(-0.2241036500, 0.4149272900)
(-0.1620185500, -0.2538387000)
(-0.1493628300, -0.0764123980)
(-0.0949611300, -0.0328289700)
(0.1539433100, -0.1519933800)
(0.1473150000, 0.5024794700)
(-0.2622958200, -0.0830622960)
(-0.1904586800, 0.1033838400)
(-0.1396389100, -0.0314963500)
(0.0522152960, 0.2171974900)
(-0.0792945060, -0.1586392600)
(-0.3143488000, -0.2731277200)
(-0.0882363920, -0.1322376700)
(-0.2273149000, 0.2255075500)
(0.0900637140, -0.2390864300)
(-0.0999437450, -0.0993648610)
(-0.0644728660, -0.1141121500)
(0.1455547100, 0.0591026700)
(-0.0310439310, 0.1077857300)
(-0.2835998400, -0.2519910400)
(-0.0385015260, 0.2455844300)
(0.0030718148, -0.1857667300)
(0.1040197700, -0.0841089700)
(0.0116337300, 0.1533349700)
(-0.2489230400, -0.2124873800)
(0.3825622000, -0.3023954700)
(0.0600763670, -0.2489919700)
(-0.1580069100, 0.3011226500)
(0.0925782590, -0.0776827260)
(-0.1447893200, 0.1721557900)
(-0.2349206700, 0.0362938380)
(-0.0566788000, -0.2740568700)
(-0.1737065500, -0.0556918060)
(-0.1394462600, 0.1657333700)
(0.1024914200, 0.1691828900)
(0.0787702200, -0.1817552200)
(0.2428397700, 0.3358801100)
(0.2555091900, -0.0308263670)
(0.3847126800, -0.2526013500)
(0.1464189000, 0.2563320400)
(-0.1156744400, -0.3136697700)
(-0.0717856400, -0.2019759200)
(-0.0688935750, -0.4544847000)
(0.0564222620, 0.4697411900)
(-0.1447681000, 0.2530908300)
(-0.1100525400, 0.0538035120)
(0.1910040500, -0.1452758700)
(-0.0928741730, 0.1807415400)
(0.1378348100, -0.1245484400)
(-0.2609359200, -0.3563058200)
(0.0525702390, -0.3382823800)
(0.1118479400, -0.1654913900)
(-0.1663148900, -0.2693724200)
(0.0204541180, 0.5428371900)
(0.0267927020, 0.5290643700)
(-0.0507842390, 0.3469720300)
(0.1614113700, 0.4197114700)
(0.4513274200, 0.0044383118)
(0.2123149400, 0.2517420100)
(0.2707771800, 0.1644441600)
(0.0011130248, 0.0855362730)
(-0.2335275800, -0.2743885400)
(0.0229873260, -0.5274120100)
(0.2766583400, -0.3868694100)
(0.2767728700, 0.0723916660)
(-0.0848099520, 0.3555661800)
(0.0868386910, -0.4571790300)
(0.0847070100, 0.0330279820)
(0.2457650100, 0.1382625500)
(0.1276670600, 0.1659382600)
(0.3348436400, 0.1971245100)
(-0.0220806830, 0.5320255900)
(0.0156731950, 0.2112138300)
(0.1134790100, -0.0262759870)
(-0.1103723300, 0.5064702200)
(0.3047502400, 0.5119987400)
(-0.2478212500, 0.0218123390)
(-0.0490311810, 0.0007404741)
(-0.0114704680, 0.0249726760)
(-0.0019821708, 0.3657865600)
(0.0017094879, 0.2413916300)
(-0.2868518600, -0.0536948130)
(0.0245211650, 0.0229647760)
(-0.0286537390, 0.2211704100)
(-0.1609470300, 0.0416872260)
(0.0857567400, -0.2500429900)
(-0.1336367200, -0.0107886600)
(-0.0878026270, -0.0786712730)
(0.0537133430, 0.0868723110)
(-0.1076369400, 0.0564932420)
(-0.0375092240, 0.4001302400)
(-0.0630157630, -0.3082400900)
(-0.0856709350, -0.1888906400)
(-0.0780491260, -0.2252342400)
(0.1868403300, 0.0708624540)
(-0.1865604600, 0.1419876400)
(-0.1640693800, 0.1711570500)
(-0.1302417800, -0.0855734840)
(0.0208864070, -0.5037741500)
(-0.2179320700, -0.0449687520)
(-0.2521392800, -0.0090909265)
(-0.2179223000, -0.0603140340)
(0.2556876400, -0.3053023000)
(-0.0746711960, 0.0631463900)
(-0.1050084900, 0.2260904800)
(0.0865260080, 0.4723587200)
(0.0406127000, 0.0913057270)
(-0.0289110610, -0.3058410100)
(0.0766304880, 0.2845338800)
(-0.4341578200, -0.3421585100)
(0.0359539000, 0.1912870100)
(-0.0150281800, -0.1599651900)
(-0.0490022530, 0.0216989740)
(0.1300506300, -0.0303066370)
(0.3120275300, 0.2469270600)
(-0.0922922400, -0.3946525800)
(0.0354435780, 0.0108225440)
(-0.3465963800, 0.0066520176)
(-0.0223007890, 0.0792623710)
(0.0584639220, 0.0137504360)
(0.1238071000, -0.1545522100)
(-0.1303809300, -0.1397448600)
(-0.0050149845, -0.1728403200)
(-0.3341818500, -0.5880271300)
(0.1193050600, -0.3602541100)
(-0.0662654680, -0.0493437560)
(-0.0439240850, 0.4257939700)
(0.0516886860, 0.0294015100)
(0.0978716950, -0.2574457100)
(-0.1156248300, -0.0993806830)
(-0.2643770400, -0.0141578030)
(-0.1775848800, -0.1622277000)
(0.1257248900, 0.1873808800)
(0.0725161290, 0.4647892600)
(0.0070266145, 0.0396909580)
(-0.2689388700, -0.4834490100)
(-0.0550347770, -0.0581880820)
(-0.2249139200, -0.2299305100)
(-0.1809029400, -0.2553419400)
(0.1969358600, 0.2717651400)
(-0.0221619220, 0.4157168300)
(-0.0889803170, 0.3989931800)
(0.0980143360, 0.2738984000)
(0.0934052020, 0.2611979700)
(-0.1610783400, -0.0646383690)
(-0.0353890720, -0.0609984690)
(0.1164233200, -0.2998246900)
(-0.0589085160, 0.4744136400)
(0.0177645420, 0.0554167000)
(-0.0673573100, -0.0747140290)
(0.3927670000, 0.2288047400)
(-0.0206500310, 0.0166509570)
(0.1295754400, 0.2333227700)
(0.0808268930, -0.0969265690)
(-0.0006537703, 0.5990355600)
(0.1037491900, 0.0201032720)
(-0.0964997240, -0.2731960600)
(-0.0107902580, 0.2254682600)
(-0.0169002840, -0.0437133670)
(0.2293280000, 0.3594157300)
(0.1001225500, 0.0305406920)
(0.0754975410, 0.1532803500)
(-0.0668426610, -0.1727375400)
(-0.0535746150, 0.1911748600)
(0.1381688400, -0.1045162700)
(-0.1445107000, 0.2475181000)
(-0.1021973200, -0.1170725700)
(-0.3089221600, -0.2357456200)
(-0.1126217500, -0.0049403234)
(-0.2142810700, -0.0038863480)
(0.0733396930, -0.4151962700)
(-0.0858293510, 0.5781205100)
(-0.1304790500, 0.0305376580)
(0.0024727513, -0.2248324300)
(-0.0454626800, -0.0040001491)
(-0.1780215500, -0.1838983500)
(-0.0514739440, -0.3698127900)
(0.1968109600, -0.2208890100)
(0.1333501700, 0.2675209400)
(0.2109687700, 0.1687502400)
(0.0910473970, 0.0483453470)
(0.0152596720, 0.0688115500)
(-0.1919755900, -0.1146773000)
(-0.0346512890, -0.1990734400)
(0.1912915800, 0.2212029200)
(0.0084386689, 0.0230193850)
(-0.0706534020, -0.2243986000)
(0.0721404130, 0.0045970127)
(-0.0888943100, 0.0482586050)
(-0.2975535700, -0.1171656100)
(-0.1849976800, 0.0391229430)
(-0.0662003880, 0.5520953800)
(0.1934905900, -0.0864485090)
(-0.0093327726, 0.1323918500)
(0.0547132220, 0.3238475700)
(0.0411420650, 0.3817891400)
(-0.0240644610, -0.0885551270)
(-0.0184489220, 0.1841217600)
(0.0780146710, 0.5675243300)
(-0.2250868000, 0.0676244760)
(-0.0887860780, -0.1122379200)
(-0.0868054470, -0.3821497000)
(-0.2290482100, -0.0491767640)
(-0.0623051280, 0.0869052380)
(-0.0397078580, 0.4902973600)
(0.2214431000, 0.2868797900)
(-0.1239105200, -0.4568448600)
(0.3067698300, -0.2673059300)
(-0.2577395500, 0.0785863180)
(0.2111712700, -0.2335912100)
(-0.0110018820, 0.0089827290)
(-0.0545089240, -0.0212171840)
(-0.0126786430, -0.3503446200)
(0.2070144800, -0.0558241740)
(-0.3053084100, -0.5823261700)
(-0.1169400000, 0.2920483700)
(-0.1474538800, -0.0772434670)
(0.2226893700, -0.3121142000)
(-0.0677441060, 0.1333680100)
(0.0402769320, 0.0236613800)
(-0.0292095870, 0.7371221100)
(0.4467931600, -0.1077330200)
(-0.0945156210, -0.1432647200)
(0.0162894150, -0.1045319400)
(-0.1408052000, -0.0069948504)
(-0.1030477600, 0.1672307200)
(0.0782366420, -0.1168280800)
(-0.1106751700, -0.2958883300)
(-0.1107816400, -0.2683028400)
(0.0244086210, -0.0721890560)
(-0.0151958450, 0.2172417200)
(0.1584319900, 0.2834893500)
(-0.3233958700, 0.1643004900)
(0.1767558500, -0.0470660590)
(0.0387423240, 0.2780336100)
(0.0943644770, -0.4654739300)
(0.2128004600, 0.2075158600)
(0.0992224830, 0.1806432000)
(0.0870709660, -0.1785395900)
(-0.1227746800, 0.1051175200)
(-0.0042323076, 0.1207397600)
(-0.0763017690, -0.0407245380)
(-0.1087349800, -0.2386057500)
(0.2265989600, -0.1263228000)
(0.1702760800, -0.0011092627)
(-0.0111712050, 0.1437278000)
(-0.1160303000, 0.1146510900)
(0.0946823610, 0.0521833760)
(0.1412400800, -0.1913159900)
(-0.0422613030, -0.3133655500)
(-0.3512770500, -0.3679726900)
(-0.1221378800, 0.3952915600)
(-0.1022108600, -0.1799596900)
(0.0464782090, 0.2815505500)
(0.2314440000, 0.0555647460)
(0.0320215900, -0.1238396700)
(0.1303230900, -0.1098699900)
(-0.1539090000, -0.2230728000)
(0.1008911700, 0.3412542100)
(0.4213916200, -0.1364930200)
(-0.3523799100, -0.1288398700)
(0.0090010154, 0.4778561000)
(0.0439557230, 0.4798511900)
(0.3593357800, 0.0676824110)
(-0.2466200900, -0.0476964010)
(0.0543364590, 0.5159439100)
(-0.1573893900, 0.3514918100)
(0.1025767600, -0.1171820900)
(0.0884577190, 0.0193814160)
(0.0852829820, -0.0535777060)
(-0.1080681300, -0.2045553700)
(-0.0747752350, -0.0425263890)
(-0.1616560600, -0.0972355260)
(0.1341420000, -0.2622426300)
(0.2917364900, -0.0215467070)
(-0.3711297000, -0.0269589490)
(-0.0026866160, -0.1888927100)
(-0.1836576300, -0.2584659000)
(-0.3350164600, -0.2090747000)
(-0.2269517500, 0.1515981400)
(-0.3085533300, -0.0482950570)
(-0.3269820300, -0.1928267300)
(-0.0000953724, -0.0421587050)
(0.0458809820, -0.0775189790)
(-0.0154968440, 0.0426714730)
(-0.1958236200, -0.2625510100)
(-0.2505143700, -0.2867206900)
(-0.1185523300, 0.4803778900)
(-0.2030390300, 0.1846290800)
(-0.1718712800, -0.1081087800)
(-0.0855573940, 0.1299218300)
(-0.1441869600, 0.0877897330)
(-0.1060965400, -0.2854970500)
(-0.2035145800, 0.0342936950)
(-0.1813461500, 0.0959690280)
(-0.0819476580, -0.2689824600)
(0.1476430400, 0.4895947300)
(0.0779098330, 0.2944872700)
(0.0036278604, 0.0893974050)
(0.0364888240, -0.5250520400)
(0.1858766800, -0.6791442800)
(0.1206678000, 0.2596681900)
(0.2335699800, -0.4287754000)
(-0.2915253300, -0.3071517800)
(-0.1733998300, 0.3978733900)
(0.0345208530, 0.3282329500)
(-0.1351481900, -0.1092032700)
(0.0387983370, -0.1674302500)
(-0.0513642230, -0.3550830500)
(-0.0812360850, 0.3520185000)
(0.0068493148, 0.4919184100)
(-0.1107010000, 0.1394738600)
(-0.0793015290, -0.1734241400)
(-0.1012216300, 0.4651525700)
(-0.0024700899, 0.2710719800)
(0.0302175530, 0.0786255630)
(-0.0673730090, -0.1165611700)
(-0.1805334100, 0.1719908200)
(0.0393620160, 0.0632827620)
(0.0843782640, -0.1855482100)
(-0.0160156380, -0.4262995600)
(-0.0783523770, -0.1224236000)
(-0.0193807530, 0.0411948850)
(-0.0009287758, -0.4284986200)
(-0.1891349800, -0.1154835000)
(0.1174518600, 0.1381641900)
(-0.0657859080, -0.2172574800)
(-0.0076180675, -0.2196837700)
(-0.0684579970, -0.2174645300)
(0.0641345650, 0.4436553600)
(-0.1500617500, 0.0964845650)
(-0.0723128950, -0.1031458100)
(-0.2003993000, -0.1017284300)
(0.1948903400, 0.4895784300)
(-0.0054811680, 0.0403979440)
(0.0056543827, 0.6338025700)
(-0.0711796440, 0.3594406800)
(0.0327253100, 0.1459739500)
(-0.1151191300, 0.3513543100)
(0.0275929530, -0.0858060670)
(-0.1565186400, 0.1168426700)
(-0.1519406600, -0.3811277400)
(-0.1008635700, 0.1965674500)
(-0.1597514100, -0.1848778000)
(-0.3152343100, -0.1999328500)
(0.1014020200, 0.3698040400)
(-0.0802394500, -0.1381820800)
(0.1211751900, 0.1034028700)
(0.1188024500, 0.3239584800)
(-0.1648441600, -0.2108809500)
(0.0191776300, -0.0931602600)
(-0.1217513200, -0.3116265700)
(-0.1131137500, -0.1094374600)
(-0.1278534500, 0.4202473600)
(0.2374191800, -0.1393434600)
(0.0501601200, 0.4705793800)
(-0.0344906080, -0.1057142100)
(-0.0137238790, 0.3294928100)
(-0.0753083540, -0.1681395600)
(-0.1688443400, 0.1472195800)
(-0.0098447225, -0.2130690700)
(0.0647697880, 0.5246150500)
(-0.0946169000, -0.4442466900)
(0.0807696510, 0.0760090190)
(-0.0946148050, 0.2155745200)
(0.1095626500, 0.4434284800)
(-0.0238265370, -0.2205543200)
(0.1373442300, -0.1136087600)
(0.1667708000, 0.0366209520)
(0.2307376900, -0.1957267100)
(0.0192567430, -0.0274300890)
};
\end{axis}
\end{tikzpicture}
}

\caption{The left plot shows the correlated errors between the elements 2 and 5 between the original tensor and the decomposed using the weighted CPD. The right plot shows the uncorrelated errors between elements 3 and 8.}
\label{fig_result1} 
\end{figure}

Starting from a uniform random $2\times2\times2$ tensor~$\Tcal$, it is decomposed using the weight matrix~\weight\ and Step~4 of Algorithm~2. This sum of rank-one tensors is denoted~$\hat{\Tcal}$. In the left plot of~\fig{fig_result1}, we plot the differences $t_5 - \hat{t}_5$ on the vertical axis, and the differences $t_2 - \hat{t}_2$ on the horizontal axis. When repeating this experiment~500 times with random tensors~$\Tcal$ and uniform random initialization points, we see that the errors are correlated, as expected.

In the right plot of~\fig{fig_result1}, the differences between $\Tcal$ and $\hat{\Tcal}$ are shown, but between the uncorrelated elements $t_3$ and~$t_8$. We see that here, the scatter plot shows uncorrelated errors between the corresponding elements of the tensor and its decomposition.

This behaviour is what is expected from the weight matrix and is discussed in the following. The weight matrix imposes extra conditions on the cost function~(\ref{eq:costfunction}) by assigning weights to the elements to be optimized. A diagonal element of \weight\ gives a weight to a single element, while an off-diagonal element of \weight\ gives a correlation between two elements. The latter implies that the errors are also correlated, which is precisely what is shown in \fig{fig_result1}.

\subsection{System identification example}
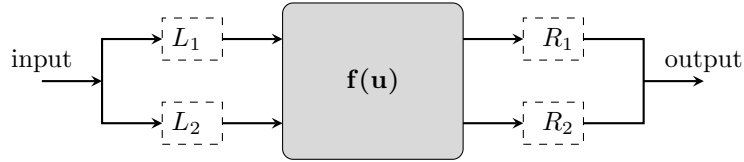
\begin{figure}[ht]
    \def\w{2.5}
    \def\ww{1.5}
    \def\wi{0.5}
    \def\fign{1}
    \def\u{\mathbf{u}}
    \def\x{\mathbf{x}}
    \def\J{\mathbf{J}}
    \def\wLN{2}
    \def\hLN{0.35}
    \centerline{\begin{tikzpicture}[>=stealth, scale=0.8]
        \filldraw[fill=gray!30, rounded corners = 4] (-\ww,-1.3) rectangle (\ww,1.3);
        \draw[dashed] (-\ww-\wLN, 0.7 - \hLN) rectangle (-\ww-1, 0.7 + \hLN);
        \draw[dashed] (-\ww-\wLN, -0.7 + \hLN) rectangle (-\ww-1, -0.7 - \hLN);
        \draw[->, thick] (-\ww-1, 0.7) -- (-\ww, 0.7);
        \draw[->, thick] (-\ww-1, -0.7)  -- (-\ww, -0.7);
        \draw[->, thick] (\ww, 0.7) -- (\ww+1, 0.7) ;
        \draw[->, thick] (\ww, -0.7) -- (\ww+1, -0.7) ;
        \draw[dashed] (\ww+\wLN, 0.7 - \hLN) rectangle (\ww+1, 0.7 + \hLN);
        \draw[dashed] (\ww+\wLN, -0.7 + \hLN) rectangle (\ww+1, -0.7 - \hLN);
        \draw[->, thick] (-\ww-1-1.5*\wLN, 0) -- (-\ww-1-\wLN, 0);
        \draw[<->, thick] (-\ww-\wLN, 0.7) -- (-\ww-1-\wLN, 0.7) --  (-\ww-1-\wLN, -0.7) -- (-\ww-\wLN, -0.7) ;
        \draw[ thick] (\ww+\wLN, 0.7) -- (\ww+1+\wLN, 0.7) --  (\ww+1+\wLN, -0.7) -- (\ww+\wLN, -0.7) ;
        \draw[->, thick] (\ww+1+\wLN, 0) -- (\ww+1+1.5*\wLN, 0);
        \draw (0,0) node {$\mathbf{f}(\u)$};
        \draw (-\ww-1-1.5*\wLN, 0) node[anchor=south] {input};
        \draw (\ww+1+1.5*\wLN, 0) node[anchor=south] {output};
        \draw (\ww+\wLN, 0.7) node[anchor = east] {$R_1$};
        \draw (\ww+\wLN, -0.7) node[anchor = east] {$R_2$};
        \draw (-\ww-\wLN, 0.7) node[anchor = west] {$L_1$};
        \draw (-\ww-\wLN, -0.7) node[anchor = west] {$L_2$};
    \end{tikzpicture}}
	\caption{Surrounding the nonlinear polynomial function \f\ by linear dynamic low-pass filters, shown by dashed lines. These operators are assumed to be known and are used when analyzing the weighted CPD on~\f. In this example, we assume the left filters $L_1$ and $L_2$ (resp. the right filters $R_1$ and $R_2$) to behave similarly and have similar transfer functions.} 
	\label{fig:sysID}
\end{figure}

In order to gain an understanding of how much impact the weight has for the CPD, the problem is examined from a distance using system identification (see \cite{Ljung1999} and \cite{Pintelon2012}). This can be done by surrounding the nonlinear function by linear low-pass filters and looking at how this system behaves with signals. \fig{fig:sysID} gives a graphical representation of the interconnection of the linear and nonlinear parts of this system. This system, with the coupled nonlinearity, comes as the result of the identification method for so-called parallel Wiener-Hammerstein system proposed in~\cite{Schoukens2015} and~\cite{SchoukensDecember10--132013}. As the name suggests, a parallel Wiener-Hammerstein system consists of parallel branches of so-called Wiener systems and Hammerstein systems. Identification methods for the latter are discussed in~\cite{Schoukens2012a} and~\cite{Schoukens2012} and for the former in \cite{Schoukens2011}.

\begin{figure}[ht]
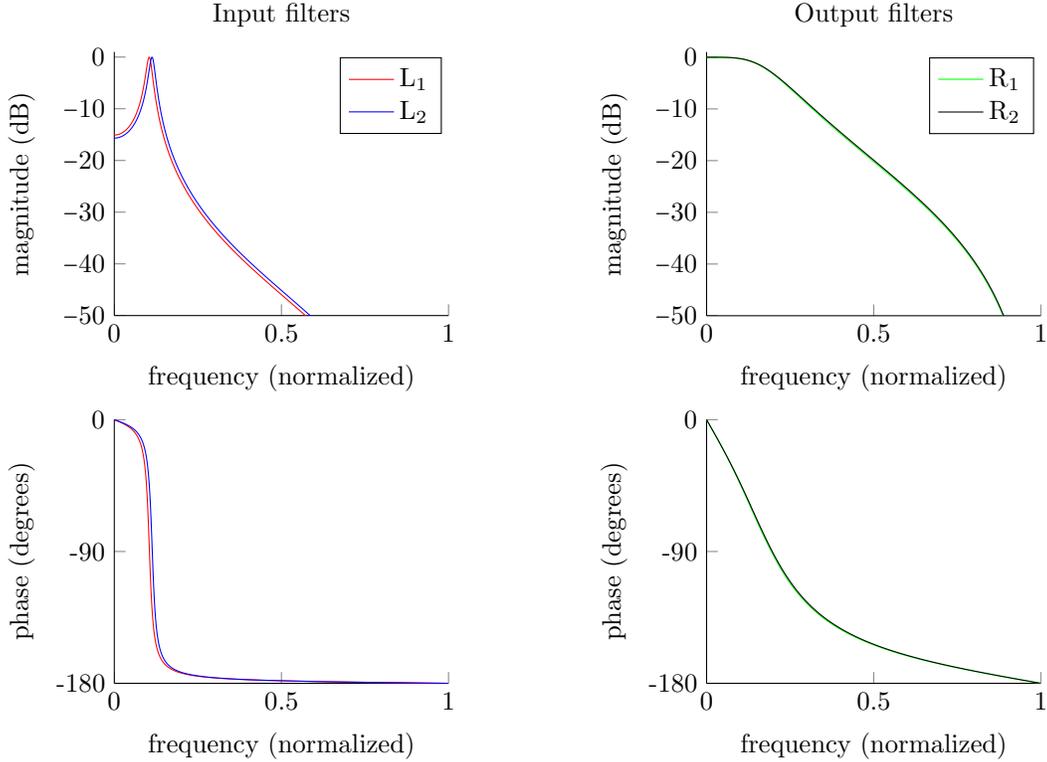

	\begin{minipage}[c]{0.5\textwidth}
	\centering
            %
%
	\end{minipage}
	\caption{The Bode plots of the input (left) and output filters (right) of the system used in this section.} 
	\label{fig_lowPass}
\end{figure}

For this example, we wish to decouple the multivariate polynomial given by 
$$
\f(\u) = \begin{cases}
0.09u_1^3 - 3.3u_1^2u_2 + 0.22u_1^2 + 5.0u_1u_2^2 - 0.44u_1u_2 - 0.25u_1 - 2.2u_2^3 + 0.41u_2^2 + 0.84u_2
\\
- 0.042u_1^3 + 3.2u_1^2u_2 - 0.21u_1^2 - 4.9u_1u_2^2 + 0.45u_1u_2 - 0.053u_1 + 2.3u_2^3 - 0.12u_2^2 - 0.27u_2
 \end{cases}.
$$
Also, we assume that the coefficients of \f\ are correlated with the covariance matrix~$\mathbf{\Sigma}_\f$ given at the top of page~\pageref{mat_covar}.
\begin{figure}
\leavevmode\kern-3em\scalebox{0.65}{
$
\mathbf{\Sigma}_\f = 
\left[
\begin{array}{@{}rrrrrrrrrrrrrrrrrr@{}}                                                                      
   \textbf{2.0} & -1.5 & 0.1 & -0.3 & 0.1 & -9.3 & 25.7 & -26.7 & 9.3 & -2.0 & 1.5 & -0.1 & 0.3 & -0.1 & 9.3 & -25.7 & 26.7 & -9.4 \\                
   -1.5 & \textbf{1.8} & -0.0 & 0.0 & -0.0 & 7.5 & -24.6 & 27.9 & -11.3 & 1.5 & -1.8 & 0.0 & -0.0 & 0.0 & -7.5 & 24.5 & -27.8 & 11.2 \\              
   0.1 & -0.0 & \textbf{9.6} & -16.0 & 7.3 & -3.7 & 6.5 & 0.3 & -2.2 & -0.1 & 0.0 & -9.6 & 15.9 & -7.3 & 3.7 & -6.6 & -0.2 & 2.2 \\                  
   -0.3 & 0.0 & -16.0 & \textbf{32.2} & -16.6 & 5.5 & -5.5 & -6.1 & 5.0 & 0.3 & -0.0 & 15.8 & -32.0 & 16.5 & -5.5 & 5.7 & 5.7 & -4.9 \\              
   0.1 & -0.0 & 7.3 & -16.6 & \textbf{10.4} & -2.7 & 1.6 & 4.2 & -2.7 & -0.1 & 0.0 & -7.3 & 16.4 & -10.3 & 2.7 & -1.8 & -4.0 & 2.6 \\                
   -9.3 & 7.5 & -3.7 & 5.5 & -2.7 & \textbf{105.7} & -275.3 & 227.3 & -56.6 & 9.3 & -7.5 & 3.7 & -5.5 & 2.7 & -105.5 & 275.4 & -228.3 & 57.2 \\      
   25.7 & -24.6 & 6.5 & -5.5 & 1.6 & -275.3 & \textbf{850.9} & -824.9 & 245.7 & -25.6 & 24.6 & -6.5 & 5.6 & -1.6 & 273.5 & -847.5 & 823.8 & -246.2 \\
   -26.7 & 27.9 & 0.3 & -6.1 & 4.2 & 227.3 & -824.9 & \textbf{937.1} & -327.3 & 26.5 & -27.9 & -0.2 & 5.9 & -4.1 & -225.3 & 819.0 & -932.4 & 326.5 \\
\end{array}\right] $
}
\label{mat_covar}
\end{figure}
In this matrix, the element on position $(i,j)$ contains the covariance between the $i$-th and $j$-th coefficient of \f. Finally, the graphical representations of the low-pass input and output filters are given in~\fig{fig_lowPass}.

\begin{figure}[ht]
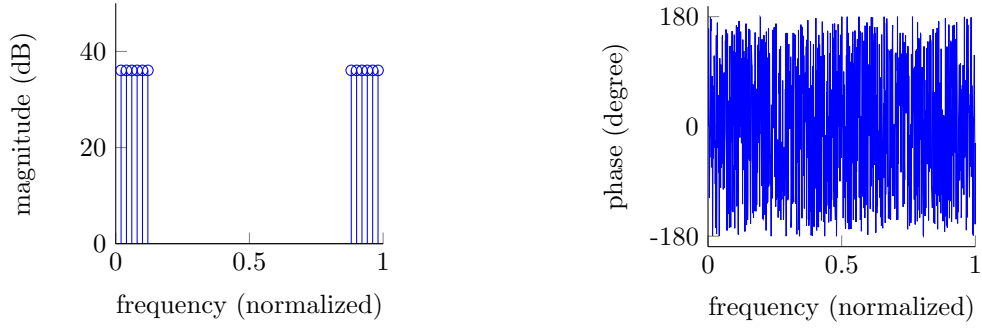

	\begin{minipage}[c]{0.5\textwidth}
	\centering
%

	\end{minipage}
	\caption{A graphical representation of the validation signal, in the frequency domain. The left plot shows the magnitude of the signal (in dB), the right plot shows the phase. The frequency axis is normalized.} 
	\label{fig_valSignal}
\end{figure}

The decoupling is done four times, in order to compare the effect of the different kinds of weights to the CPD:
\begin{enumerate}
\item Decoupling without any weight, using the method of~\cite{Dreesen2015},
\item Decoupling with the element-wise weight,
\item Decoupling with the slice-wise weight,
\item Decoupling with the dense weight.
\end{enumerate}
In order to compare the results, a validation signal is sent through the original coupled system of~\fig{fig:sysID} and also through the four different decouplings. This validation signal is a random-phase multisine (see \cite{Pintelon2012}), as shown in~\fig{fig_valSignal}.

\begin{figure}[ht]
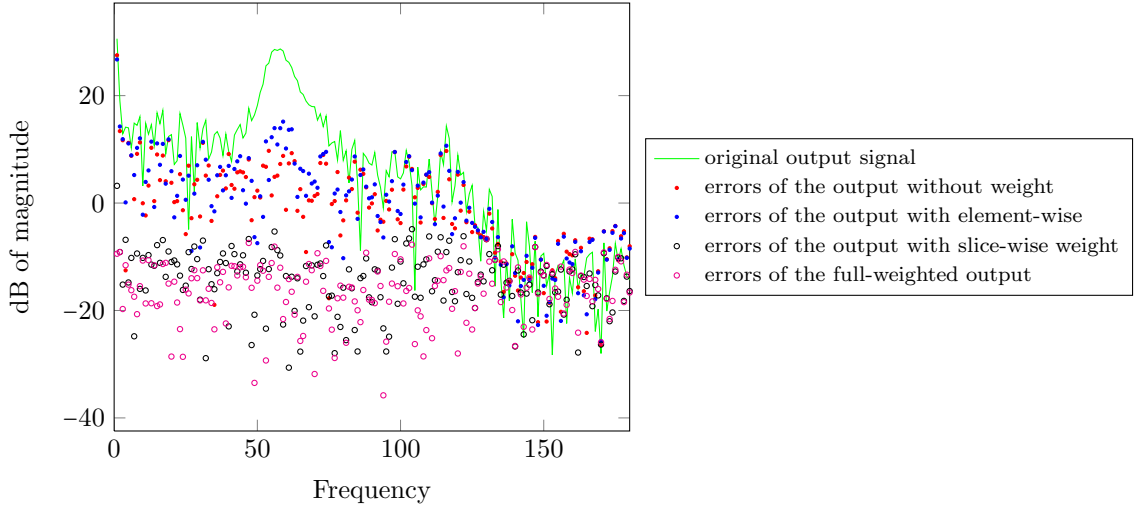

	\begin{center}

        \end{center}
	\caption{The errors between the the original output signal and the output by the different decoupling methods. Plotted in the frequency domain; the vertical axis is in dB of the magnitude of the signals.} 
	\label{fig:results}
\end{figure}

The output signal and the errors between the different decouplings are plotted in~\fig{fig:results}. This plot shows the magnitude of the output and output errors, in dB, with respect to the frequency of the signals. From this plot, it is clear that a slice-wise or dense weight reduces the errors significantly with respect to no weight or element-wise weight.

Next to these numerical experiments, extensive simulations were done with multiple different systems and several coupled multivariate polynomials. From these, general observations can be made: the slice-wise and dense weight decompositions are at least as good as the decomposition with no weight described in~\cite{Dreesen2015}. Overall, the element-wise weight does not incorporate enough information to improve~\cite{Dreesen2015} and has comparable results. Finally, for difficult decomposition problems where~\cite{Dreesen2015} does not work well, we observe improvements with the slice-wise or dense weighted decomposition.

\section{Conclusion and future work}\label{sec:conclusion}
In this paper, our starting point is a coupled representation of a multivariate polynomial~$\f$, which does not have an exact decomposition with the prespecified number of branches. It is assumed, however, that the covariance matrix~$\mathbf{\mathbf{\Sigma}}_\f$ of the coefficients of~$\f$ is known before the decoupling process. We have then generalized the decoupling algorithm described in~\cite{Dreesen2015} to this noisy case, by considering a weight factor in the Canonical Polyadic Decomposition. Three weight factors have been considered, based on three different covariance matrices: (1) an element-wise covariance matrix, (2) a slice-wise covariance matrix and (3) a dense covariance matrix. In cases (1) and (2), the matrices are of full rank, while the matrix in case (3) is rank-deficient. That is why extra equations are found for the third case, based on a SVD of the dense covariance matrix.

The results are promising and at least as good as the unweighted decoupling method described in~\cite{Dreesen2015}. When considering the decoupling problem inside the framework of linear filters, improvements are observed and discussed.

As future work, we would like to investigate how to make approximations of polynomials with a high number of branches by polynomials with a lower number of branches. Also, generalizations to other basis functions can be studied, and in what way these modify the proposed method.

\appendix
\section{Derivation of~(\ref{eq_lowRankALS})}\label{sec:app1}
In this section, we show the derivations for finding the solution~(\ref{eq_lowRankALS}). This will be done in a more general case, where we consider the following weighted least squares problem:
\begin{equation}
\min_\x \norm{\y - \A \x}_{\weight}^2.
\label{eq_wlq}\end{equation}
The weight matrix \weight\ is defined as the pseudo-inverse of the rank-deficient covariance matrix $\mathbf{\Sigma}\in\ER^{mnN\times mnN}$ of rank $\rSigma$.
Let $\D_1$ and $\U_1$ be defined as in Section~\ref{sec:fWeight_part1}, starting from the singular value decomposition of $\mathbf{\Sigma} = \U\D\U^T$:
\begin{equation}
\rule[0pt]{0pt}{23pt}
\U = \Bigl[
\begin{array}{@{}c|c@{}}
		 \smash{\overbrace{\U^{(1)}}^{\makebox[0pt]{$\scriptstyle mnN\times \rSigma$}}}
 &
 		 \smash{\underbrace{\U^{(2)}}_{\makebox[0pt]{$\scriptstyle mnN\times(mnN-\rSigma)$}}}
\end{array}
\Bigr]
\quadtext{and} \D = \left[\begin{array}{@{}c|c@{}}
\smash{\overbrace{\D^{(1)}}^{\rSigma\times\rSigma}}
 & 0 \\
\hline
0 & 0
\end{array} \right].
\label{eq_V1V2}\end{equation}
We then define
$$
\Q = (\sqrt{\D^{(1)}})^{-1} (\U^{(1)})^T,
$$
such that $\Q^T\Q = \weight$.

We note that minimizing expression~(\ref{eq_wlq}) can be rewritten as 
$$
	\min (\e^T \,\weight \,\e),
$$
where $\e = \A \x - \y$ is the error between modeled output and real output. If we denote the transformed error by $\tilde{\e} = \Q \, \e = \Q \, \A \, \x - \Q \,\y$, then this minimization becomes
$$
\e^T \,\weight \, \e = \e^T \Q^T \Q \, \e = \tilde{\e}^T \tilde{\e}.
$$
This implies that the weighted least squares solution from (\ref{eq_wlq}) can be transformed as the (unweighted) least squares solution of $\tilde{\e}^T \tilde{\e}$. This is given by 
$$
\hat{\x} = (\Q \A)^{\!\dagger} (\Q y).
$$
This is precisely the form used in the expression~(\ref{eq_lowRankALS}) and is what we wanted to prove.

\section{Derivation of~(\ref{eq_resultEquation2})}\label{sec:app2}
The equations of Appendix~\ref{sec:app1} depend on the $\U^{(1)}$ part of the matrix $\U$. Here, we will find extra equations using~$\U^{(2)}$. For this, we assume that the noise $\v$ added to the model
\begin{equation}
\y = \A \x + \v
\label{eq_model}\end{equation}
is correlated, such that $\v = \T \v_\text{uncor}$. Here, $\v_\text{uncor}$ is independent, identically distributed noise with $\var(\v_\text{uncor}) = \sigma^2$ and $\T\in\ER^{N\times M}$ creates correlations between elements of $\v_\text{uncor}$. Finally, we assume that $N>M$.

With these notations, we compute the covariance matrix $\mathbf{\Sigma}$ of $\v$
$$
\mathbf{\Sigma} = \cov(\v) = \E[\v \,\v^T] = \sigma^2 \,\T\,\T^T.
$$
Because the rank of $\mathbf{\Sigma}$ is at most $M$, then $\mathbf{\Sigma}$ is rank-deficient. Using the same notations as in the equations~(\ref{eq_V1V2}), it follows that
\begin{equation}
(\U^{(2)})^T \v = (\U^{(2)})^T  \, \T \, \v_\text{uncor} = 0.
\label{eq_V2}\end{equation}
This follows from the that fact $\U^{(2)} \perp \U^{(1)}$ and thus that $\U^{(2)} \perp \T$. Finally, substituting Equation~(\ref{eq_model}) in Equation~(\ref{eq_V2}) gives 
$$
(\U^{(2)})^T \y - (\U^{(2)})^T \A \x = 0,
$$
which has as approximated solution for $\x$:
$$
\hat{\x} = \bigl((\U^{(2)})^T \A\bigr)^\dagger \, \bigl((\U^{(2)})^T \y\bigr).
$$
This approximation is used for finding extra equations in Section~\ref{sec:fWeight_part2}.


\end{document}